\documentclass[11pt,a4paper]{article}
\usepackage[top=1.0in, bottom=1.0in, left=1.0in, right=1.0in]{geometry}
\usepackage{setspace}
\usepackage{amsfonts}
\usepackage{amssymb}
\usepackage{graphicx}
\usepackage{amsfonts}
\usepackage{fancyhdr}
\usepackage{float}
\usepackage{subfigure}
\usepackage{epstopdf}
\usepackage{authblk}
\usepackage[T1]{fontenc}
\usepackage{algorithm}
\usepackage{algorithmic}
\usepackage{amsthm, amsmath, enumitem}
\usepackage{color}
\usepackage{multirow}
\setlist{nosep, topsep=2pt, partopsep=0pt, parsep=0pt, itemsep=0pt}

\usepackage[T1]{fontenc}
\usepackage{authblk}
\usepackage[colorlinks,linktocpage,linkcolor=blue]{hyperref}
\usepackage{amsfonts}
\usepackage{amssymb}
\usepackage{graphicx}
\usepackage{amsfonts}
\usepackage{fancyhdr}
\usepackage{float}
\usepackage{subfigure}
\usepackage{epstopdf}
\usepackage{algorithm}
\usepackage{algorithmic}
\usepackage{amsthm, amsmath, enumitem}
\usepackage{color}
\usepackage{multirow}

\newtheorem{theorem}{Theorem}[section]
\newtheorem{lemma}{Lemma}[section]

\newtheorem{Algo}{Algorithm}[section]
\newtheorem{example}{Example}

\def\bq{\begin{equation}}
\def\eq{\end{equation}}
\def\br{\begin{eqnarray}}
\def\er{\end{eqnarray}}
\def\brr{\bq\begin{array}{rlll}}
\def\err{\end{array}\eq}

\def\text#1{\hbox{#1}}

\newcommand{\bsub}{\begin{subequations}}
\newcommand{\esub}{\end{subequations}$\!$}

\graphicspath{{Figure3/}} 

\usepackage{bm}
\usepackage{mathrsfs}%
\usepackage{cite}
\numberwithin{equation}{section}
\begin{document}
\title{\Large{Solving the inverse Source Problems for wave equation with final time measurements by a data driven approach}}

\author[1]{Qiling Gu\thanks{guql@sustech.edu.cn}}
\author[1]{Wenlong Zhang\thanks{zhangwl@sustech.edu.cn}}
\author[2,3]{Zhidong Zhang\thanks{zhangzhidong@mail.sysu.edu.cn}}
\affil[1]{\normalsize{Department of Mathematics, Southern University of Science and Technology (SUSTech), Shenzhen, Guangdong, China.}}
\affil[2]{\normalsize{School of Mathematics (Zhuhai), Sun Yat-sen University, Zhuhai 519082, Guangdong, China}}
\affil[3]{\normalsize{Guangdong Province Key Laboratory of Computational Science, Sun Yat-sen
University, Guangzhou 510000, Guangdong, China}}
\maketitle

\begin{abstract}
This paper develops a discrete data-driven approach for solving the inverse source problem of the wave equation with final time measurements. Focusing on the $L^2$-Tikhonov regularization method, we analyze its convergence under two different noise models, using noisy discrete spatial observations. By exploiting the spectral decomposition of the forward operator and introducing a noise separation technique into the variational framework, we establish error bounds for the reconstructed solution $u$ and the source term $f$ without requiring classical source conditions. Moreover, an expected convergence rate for the source error is derived in a weaker topology. We also extend the analysis to the fully discrete case with finite element discretization, showing that the overall error depends only on the noise level, regularization parameter, time step size, and spatial mesh size. These estimates provide a basis for selecting the optimal regularization parameter in a data-driven manner, without a priori information. Numerical experiments validate the theoretical results and demonstrate the efficiency of the proposed algorithm.
\end{abstract}

\date{}

\hskip\parindent
\section{Introduction}\label{secInt}

 Wave equations serve as fundamental models for wave propagation phenomena across geophysical exploration, non-destructive testing, and medical imaging. The accurate reconstruction of source terms from boundary or internal measurements is essential for identifying hidden physical structures \cite{16,17,18,27,28,33,43,44,45,46}. This paper develops a discrete data-driven approach for solving the inverse source problem of the wave equation with final-time measurements. We establish rigorous convergence estimates for the $L^2$-Tikhonov regularized solution, deriving error bounds for both the continuous setting and its finite element discretization, using scattered pointwise data contaminated by noise.

Recent methodological advances in wave-induced inverse source problems include quasi-reversibility schemes for bounded cavities \cite{44}, multi-frequency factorization techniques \cite{43}, moving source localization methods \cite{17}, hybrid identification approaches \cite{45}, and theoretical analyses for biharmonic waves \cite{46}. A common feature among existing convergence analyses, such as those in \cite{46}, is their reliance on source conditions—typically Sobolev smoothness assumptions—to establish theoretical guarantees. In contrast, this work introduces a new analysis framework that achieves convergence without imposing any source conditions. By leveraging a discrete data-driven setup with noisy final-time measurements, we derive error estimates in weaker topologies, effectively bypassing the limitations of condition-dependent theories.

A key contribution of this work is the derivation of finite element (FE) error estimates for the inverse source problem—a topic that remains largely unexplored for wave equations due to inherent challenges in hyperbolic system analysis. While some error estimates exist for inverse elliptic and parabolic problems \cite{2,30}, analogous results for wave equations are scarce. We present the first FE error analysis for the standard $L^2$
-Tikhonov regularization of the wave equation inverse source problem. Our approach accommodates realistic noise models and provides explicit error dependencies on four computational parameters: the noise level $\sigma$, the regularization parameter $\alpha$, the spatial mesh size $h$, and the temporal step size $\tau$. These estimates support a fully data-driven parameter selection strategy, eliminating the need for a priori information.

Methodologically, we adopt a scattered data framework \cite{1,2,3,4,5,6,7,8} in which measurements are taken at discrete points without additional smoothness assumptions. Using a novel noise separation technique within the variational formulation, we extend the scattered data analysis introduced in \cite{29}—originally formulated for deterministic settings—to accommodate noisy measurement environments. Under standard compatibility conditions on the temporal source term $g(t)$, we establish regularity estimates for the forward solution, spectral properties of the forward operator, and uniqueness results for the inverse problem.

The main contributions of this paper are:
\begin{enumerate}
\item[(i)] Derivation of discrete finite element error estimates without source conditions, with explicit dependence on the noise level, regularization parameter, temporal step size, and spatial mesh size;
\item[(ii)] Development of a self-consistent, data-driven algorithm for optimal regularization parameter selection without prior knowledge.
\end{enumerate}

The remainder of the paper is organized as follows. Section~\ref{1} introduces preliminaries, including regularity estimates for the forward solution, spectral properties of the forward operator, and empirical error bounds from finite element discretization. Section~\ref{2} analyzes convergence of the continuous regularized solution under two noise models. Section~\ref{3} extends the analysis to the fully discrete case. Section~\ref{4} presents a self-consistent algorithm for choosing $\alpha$, accompanied by numerical experiments.

Let $W^{m,p}(\Omega)$ denote the Sobolev space for integer $m \geq 0$ and $p \geq 1$, with norm $\|\cdot\|_{W^{m,p}(\Omega)}$; when $p = 2$, we write $H^{m}(\Omega)$ with norm $\|\cdot\|_{H^{m}(\Omega)}$. The space of continuous functions on $\overline{\Omega}$ is denoted by $C(\overline{\Omega})$. The discrete measurement functional $\mathcal{M}$ for $v \in C(\overline{\Omega})$ is defined as:
\[
(\mathcal{M},v)_{n} := \frac{1}{n}\sum_{i=1}^{n}\mathcal{M}_{i}v(x_{i}),
\]
and the empirical inner product and norm are given by:
\[
(u,v)_{n} := \frac{1}{n}\sum_{i=1}^{n}u(x_{i})v(x_{i}),\quad \|u\|_{n} := \left(\frac{1}{n}\sum_{i=1}^{n}u^{2}(x_{i})\right)^{1/2}.
\]
Throughout, $C$ denotes a generic constant independent of $h$, $n$, and $\alpha$.

\section{Preliminaries}\label{1}
This section establishes the mathematical setting for the inverse source problem. We define the forward model, the measurement process, and the regularization strategy, followed by key results on the forward operator's properties and its numerical discretization.

\subsection{Problem Formulation and Mathematical Setting}

Consider the wave equation:
\begin{equation}\label{2.1}
\begin{cases}
u_{tt} - \Delta u = f(x)g(t), & (x,t) \in \Omega \times (0, T), \\
u(x,t) = 0, & (x,t) \in \partial\Omega \times (0, T), \\
u(x,0) = 0, \quad u_t(x,0) = 0, & x \in \Omega,
\end{cases}
\end{equation}
where \(\Omega \subset \mathbb{R}^d\) ($d=1,2,3$) is a bounded domain with a $C^2$ boundary or convex polyhedral structure. The forward problem is to compute the wave field $u(\cdot, t)$ given the source $f(x)$. The inverse problem is to recover an unknown $f \in L^2(\Omega)$ from final-time observations $u(\cdot, T)$.

For $f \in L^2(\Omega)$ and $g \in H^4(0,T)$ with $g^{(k)}(0) = 0$ for $k = 0,1,2,3$, the solution regularity is $u \in H^2(\Omega)$ (see Lemma \ref{L5.1}). This allows us to define the forward operator $G: L^2(\Omega) \to H^2(\Omega)$ as $Gf := u(\cdot, T)$.

Data are collected at discrete sensor locations $\{x_i\}_{i=1}^n \subset \Omega$ ($x_i \neq x_j$ for $i \neq j$) that are quasi-uniformly distributed, meaning there exists $B > 0$ such that
\begin{equation}\label{2.6b}
d_{\max} := \sup_{x \in \Omega} \min_{1 \leq i \leq n} |x - x_i|, \quad
d_{\min} := \min_{1 \leq i \neq j \leq n} |x_i - x_j|, \quad
\frac{d_{\max}}{d_{\min}} \leq B.
\end{equation}
The measurement model with additive noise is
\begin{equation}\label{2.3}
m_i = Gf^*(x_i) + e_i, \quad i = 1, \dots, n,
\end{equation}
where $f^*$ is the true source, $m = (m_1, \dots, m_n)^T$ is the observation vector, and $e = (e_1, \dots, e_n)^T$ is a noise vector. We consider two noise models:
\begin{itemize}
    \item[(Y1)] Independent noise with $\mathbb{E}[e_i] = 0$ and $\mathbb{E}[e_i^2] \leq \sigma^2$.
    \item[(Y2)] Independent sub-Gaussian noise.
\end{itemize}

To recover $f^*$, we use $L^2$-Tikhonov regularization. The estimator $f_n \in L^2(\Omega)$ minimizes:
\begin{equation}\label{2.2}
\min_{f \in L^2(\Omega)} \|Gf - m\|_n^2 + \alpha \|f\|_{L^2(\Omega)}^2,
\end{equation}
where $\alpha > 0$ is the regularization parameter.

\subsection{Regularity and Spectral Property of the Forward problem}
The  convergence of the Tikhonov solution is intimately tied to the analytical properties of the forward operator $G$. First, we establish the solution regularity, which underpins the subsequent analysis.
\begin{lemma}\label{L5.1}
Let $u$ be the solution to \eqref{2.1}. Suppose $f \in L^2(\Omega)$ and $g \in H^4(0,T)$ with $g^{(k)}(0) = 0$ for $k=0,1,2,3$. Then $u$ satisfies the estimate
\[
\|u\|_{L^\infty(0,T; H^2(\Omega))} + \|u_t\|_{L^2(0,T; H^2(\Omega))} + \|u_{ttt}\|_{L^2(0,T; L^2(\Omega))} + \|u_{tttt}\|_{L^2(0,T; L^2(\Omega))} \le C \|f\|_{L^2(\Omega)},
\]
where $C > 0$ depends on $\Omega$, $T$, and $\|g\|_{H^4(0,T)}$.
\end{lemma}

\begin{proof}
 Let $v(t) = (u(\cdot, t), u_t(\cdot, t))^T$. The problem  \eqref{2.1} can be written as:
\[
\begin{cases}
    v_t = A v + F(t), & \text{in } \Omega \\
    v(0) = 0,
\end{cases}
\]
where we operate in the Hilbert space $X = H_0^1(\Omega) \times L^2(\Omega)$ and define the operator $A$ and the forcing term $F(t)$ as:
\[
A = \begin{pmatrix} 0 & I \\ \Delta & 0 \end{pmatrix}, \quad F(t) = \begin{pmatrix} 0 \\ f(\cdot) g(t) \end{pmatrix}.
\]
The operator $A$ with domain $D(A) = (H^2(\Omega) \cap H_0^1(\Omega)) \times H_0^1(\Omega)$ is the generator of a $C_0$-semigroup of contractions $e^{tA}$ on $X$. The solution is given by the variation of constants formula:
\[
v(t) = \int_{0}^{t} e^{(t-s)A} F(s) \, ds.
\]
Utilizing integration by parts and the hypothesis and the hypothesis $g^{(k)}(0) = 0$ for $k=0,1,2,3$, we obtain:
\[
v^{(k)}(t) = \int_{0}^{t} e^{(t-s)A} F^{(k)}(s) \, ds, \quad \text{for } k = 1, 2, 3, 4.
\]
From the differential equation $v_t = Av+F$, we have
\[
\|\Delta u(\cdot, t)\|_{L^2(\Omega)} \leq \|Av(t)\|_X \leq \|v_t(t)\|_X + \|F(t)\|_X,
\]
which implies the first bound:
\[
\|u\|_{L^\infty(0,T; H^2(\Omega))} \le C \|f\|_{L^2(\Omega)}.
\]
An analogous argument for 
to $Av_t = v_{tt} - F_t$ yields  an $L^\infty(0,T; H^2(\Omega))$ estimate for $u_t$, implying:
\[
\|u_t\|_{L^2(0,T; H^2(\Omega))} \le C \|f\|_{L^2(\Omega)}.
\]

From the definition $v_{ttt}(t) = (u_{ttt}, u_{tttt})^T$, the first component gives the bound $\|u_{ttt}(\cdot,t)\|_{L^2(\Omega)} \le C \|u_{ttt}(\cdot,t)\|_{H_0^1(\Omega)} \le C\|v_{ttt}(t)\|_X$ via Poincaré's inequality. For its $L^2(0,T; X)$
-norm:
\[
\|v_{ttt}\|_{L^2(0,T; X)} \le C \|F^{(3)}\|_{L^2(0,T; X)} = C \|f g^{(3)}\|_{L^2(0,T; L^2(\Omega))} = C \|f\|_{L^2(\Omega)} \|g^{(3)}\|_{L^2(0,T)}.
\]
This directly provides the bound for $\|u_{ttt}\|_{L^2(0,T; L^2(\Omega))}$.

Similarly, for $\|u_{tttt}\|_{L^2(0,T; L^2(\Omega))}$, a parallel argument using $v_{tttt}(t) = (u_{tttt}, u_{ttttt})^T$ yields:

\[
\|
u_{tttt}\|_{L^2(0,T; L^2(\Omega))} \leq C\|f\|_{L^2(\Omega)}\|g^{(4)}\|_{L^2(0
,T)}.
\] 
Combining all the derived bounds completes the proof.
\end{proof}

With this regularity established, we can analyze the spectral properties of the forward operator $G$, which are key to understanding the inverse problem's structure.
\begin{lemma}\label{T3.1}
 Let $\Omega \subset \mathbb{R}^d$ be a bounded smooth domain, and let $g \in H^4(0,T)$ satisfy $g^{(j)}(0)=0$ for $j=0,1,2,3$, $g'(t)>0$, and $T>0$. Then the null space $N(G)=\{0\}$, and the eigenvalue problem
\[
(\psi, v)_{L^2(\Omega)} = \eta (G\psi, Gv), \quad \forall v \in L^2(\Omega)
\]
admits a countable set of positive eigenvalues $\eta_1 \le \eta_2 \le \dots \to \infty$ for which there exist constants $C, C' > 0$ such that
\[
C k^{4/d} \leq \eta_k \leq C' k^{4/d}, \quad \forall k \in \mathbb{N}.
\]
\end{lemma}

\begin{proof}
Let $\{\phi_n\}$ be the orthonormal eigenbasis of $-\Delta$ in $L^2(\Omega)$ with eigenvalues $0 < \mu_1 \le \mu_2 \le \cdots$. Expanding $u(x,t) = \sum u_n(t)\varphi_n(x)$ and $f(x) = \sum f_n\varphi_n(x)$ with $f_n = (f,\varphi_n)$, we obtain for each coefficient:
\begin{equation}\label{eq:ode}
u_n''(t) + \mu_n u_n(t) = f_n g(t), \quad u_n(0) = u_n'(0) = 0.
\end{equation}
By Duhamel's principle, the solution is given by
\begin{equation}\label{2.5a}
u_n(t) = \frac{f_n}{\sqrt{\mu_n}} \int_0^t \sin(\sqrt{\mu_n}(t-s)) g(s) \, ds.
\end{equation}
Assume $f \in L^2(\Omega)$ satisfies $Gf = u(\cdot, T) = 0$. Then $u_k(T) = f_k \alpha_k = 0$ for all $k$, where
\[
\alpha_k = \frac{1}{\sqrt{\mu_k}} \int_0^T \sin(\sqrt{\mu_k}(T - s)) g(s)\, ds.
\]
 Let $t_i = i\Pi_h/\sqrt{\mu_k}$ for $i=0, \dots, n$ partition the interval $(0,T]$, with $t_{n+1}=T$. Since $g(t) \in C[0,T]$ and $\sin(\sqrt{\mu_k}s)$ is sign-preserving on each subinterval $(t_i, t_{i+1}]$, the Mean Value Theorem for integrals implies there exists $\xi_i \in (t_i, t_{i+1}]$ such that
\begin{align*}
\sqrt{\mu_k}\alpha_k &= \int_0^T g(T-s)\sin(\sqrt{\mu_k}s)ds = \sum_{i=0}^n g(T-\xi_i) \int_{t_i}^{t_{i+1}} \sin(\sqrt{\mu_k}s) ds \\
  &= \sum_{i=0}^{n-1} g(T-\xi_i) \frac{2 \cdot (-1)^i}{\sqrt{\mu_k}} + g(T-\xi_n) \left( \frac{(-1)^n - \cos(\sqrt{\mu_k}T)}{\sqrt{\mu_k}} \right).
\end{align*}
If $n$ is even:
\begin{align*}
\sqrt{\mu_k}\alpha_k &= \sum_{m=0}^{\frac{n-2}{2}} \left[ g(T-\xi_{2m}) - g(T-\xi_{2m+1}) \right] \frac{2}{\sqrt{\mu_k}} 
   + g(T-\xi_n) \left( \frac{1 - \cos(\sqrt{\mu_k}T)}{\sqrt{\mu_k}} \right).
\end{align*}
If $n$ is odd:
\begin{align*}
\sqrt{\mu_k}\alpha_k &= \sum_{m=0}^{\frac{n-3}{2}} \left[ g(T-\xi_{2m}) - g(T-\xi_{2m+1}) \right] \frac{2}{\sqrt{\mu_k}} \\
  &\quad + \left(g(T-\xi_{n-1}) \frac{2}{\sqrt{\mu_k}} + g(T-\xi_n) \left( \frac{-1 - \cos(\sqrt{\mu_k}T)}{\sqrt{\mu_k}} \right) \right).
\end{align*}
Since $g'(t) > 0$ and $T>0$, we have $\alpha_k > 0$ for all $k$. Then $u_k(T) = f_k \alpha_k = 0$ implies $f_k=0$ for all $k$, so $f=0$ and $N(G)=\{0\}$.

Standard energy estimates show that $G: L^2(\Omega) \to H^2(\Omega)$ is bounded. By the Rellich-Kondrachov theorem, the embedding $H^2(\Omega) \hookrightarrow L^2(\Omega)$ is compact, so $G: L^2(\Omega) \to L^2(\Omega)$ is compact. Since $N(G)=\{0\}$, $G^*G$ is compact, self-adjoint, and strictly positive definite. The Spectral Theorem provides an orthonormal basis $\{\psi_k\}$ and eigenvalues $\beta_k \downarrow 0$ such that 
\begin{equation}\label{2.7ab}
(G^*G)\psi_k = \beta_k\psi_k\quad \text{in } L^2(\Omega), \quad (\psi_k, \psi_l)_{L^2(\Omega)}= \delta_{kl}, \quad (G\psi_k, G\psi_l) = \beta_k \delta_{kl}, \quad k,l = 1, 2, \dots,
\end{equation}

The eigenvalue problem $(\psi, v) = \eta (G\psi, Gv)$ is equivalent to $(G^*G)\psi = \eta^{-1}\psi$, so $\eta_k = \beta_k^{-1}$. Since $G^*G$ and $-\Delta$ commute, they share eigenfunctions, so $\psi_k = \varphi_k$. Thus $\beta_k = \|G\varphi_k\|^2 = \alpha_k^2$ and $\eta_k = \alpha_k^{-2}$. Applying integration by parts to $\alpha_k$ and using $g^{(j)}(0) = 0$ for $j=0,1,2,3$ yields:

\begin{equation}\label{2.5abc}
\alpha_k = \frac{g(T)}{\mu_k} - \frac{g''(T)}{\mu_k^2} + \frac{1}{\mu_k^{5/2}} \int_0^T \sin(\sqrt{\mu_k}(T - \tau)) g^{(4)}(\tau) d\tau.
\end{equation}

This expansion shows $\alpha_k \le C_2/\mu_k$ for some $C_2 > 0$. Since $g'(t)>0$ and $g(0)=0$, we have $g(T)>0$. As $k \to \infty$, $\mu_k \to \infty$ and $\mu_k \alpha_k \to g(T)$. Since $\alpha_k > 0$ for all $k$, there exists $C_1 > 0$ such that $\mu_k \alpha_k \ge C_1$ for all $k$, so $\alpha_k \ge C_1/\mu_k$.

Combining the bounds on $\eta_k = \alpha_k^{-2}$ with Weyl's law ($\mu_k = C_W k^{2/d}$) completes the proof.
\end{proof}

\subsection{Finite Element Discretization and Empirical Error}

For numerical implementation, we discretize \eqref{2.1} using conforming linear finite elements on a quasi-uniform mesh $\mathcal{M}_h$, with space $V_h \subset H_0^1(\Omega)$, and a second-order temporal scheme on a uniform time grid $t^i = i\tau$ where $\tau=T/N$. The fully discrete scheme seeks $u_h^i \in V_h$ for $i=1,\dots,N$ such that for all $v_h \in V_h$:
\begin{equation}\label{5.2}
\begin{aligned}
\frac{1}{\tau}(q_h^i - q_h^{i - 1}, v_h) + a\left(\frac{u_h^i + u_h^{i - 1}}{2}, v_h\right) &= \left(\frac{f^i + f^{i - 1}}{2}, v_h\right), \\
\frac{u_h^i - u_h^{i - 1}}{\tau} &= \frac{q_h^i + q_h^{i - 1}}{2},
\end{aligned}
\end{equation}
where $a(u_h, v_h) := (\nabla u_h, \nabla v_h)$. This scheme produces a numerical approximation to the wave field at the final time.

We define the numerical forward solution as $G_{\tau,h} f := u_h^N$. The estimator $f_h \in V_h$ is the minimizer of the fully discrete Tikhonov functional:
\begin{equation}\label{5.6}
\min_{f \in V_h} \| G_{\tau,h} f - m \|_n^2 + \alpha \| f \|_{L^2(\Omega)}^2.
\end{equation}
This $f_h$ satisfies the weak formulation:
\begin{equation}\label{2.11a}
\alpha (f_h, v_h) + (G_{\tau,h} f_h, G_{\tau,h} v_h)_n = (m, G_{\tau,h} v_h)_n \quad \forall v_h \in V_h.
\end{equation}

Finally, to connect our analysis of the continuous and discrete problems, we quantify the error between the true forward operator $G$ and its numerical approximation $G_{\tau,h}$.

\begin{lemma}\label{L2.3a}
Let $f \in L^2(\Omega)$, and suppose $g \in H^4([0,T])$ with $g^{(k)}(0) = 0$ for $k=0,1,2,3$. Then the fully discrete solution $G_{\tau,h} f = u_h^N$ satisfies
\[
\| Gf - G_{\tau,h} f \|_n \le C(h^2 + \tau^2) \| f \|_{L^2(\Omega)},
\]
where $C > 0$ is independent of $h, \tau$, and $f$.
\end{lemma}

\begin{proof}
Let $\Pi_h: C(\bar{\Omega}) \to V_h$ be the linear finite element interpolation operator. By the triangle inequality and standard interpolation estimates\cite{32} \begin{align*}
\|Gf - \Pi_h(Gf)\|_{L^\infty(K)} &\le Ch^{2-d/2}\|Gf\|_{H^2(K)} \quad \forall K \in \mathcal{M}_h, \\
\|Gf - \Pi_h(Gf)\|_{L^2(K)} &\le Ch^2\|Gf\|_{H^2(K)} \quad \forall K \in \mathcal{M}_h,
\end{align*}Let $T_K = \{x_i : x_i \in K, 1 \le i \le n\}$, we have
\begin{align*}
\| Gf - \Pi_h (Gf) \|_n^2 &= \frac{1}{n} \sum_{K \in \mathcal{M}_h} \sum_{x_i \in \mathbb{T}_{K}} |Gf(x_i) - \Pi_h Gf(x_i)|^2 \\
&\le \frac{1}{n} \sum_{K\in \mathcal{M}_h} \# \mathbb{T}_{K} \| Gf - \Pi_h Gf \|_{L^\infty(K)}^2 \le C h^4 \| Gf \|_{H^2(\Omega)}^2,
\end{align*}
where we used the quasi-uniformity of the mesh $\mathcal{M}_h$ and sensors $\{x_i\}_{i=1}^n$, which implies $\# \mathbb{T}_{K} \le C n h^d$.
Using local inverse estimates, we can also bound the numerical solution's interpolation error:
\begin{align*}
\|G_{\tau,h}f - \Pi_h(Gf)\|_n^2 &\le \frac{1}{n} \sum_{K \in \mathcal{M}_h} \#T_K \|G_{\tau,h}f - \Pi_h(Gf)\|_{L^\infty(K)}^2 \\
&\le \frac{1}{n} \sum_{K \in \mathcal{M}_h} \#T_K |K|^{-1} \|G_{\tau,h}f - \Pi_h(Gf)\|_{L^2(K)}^2 \\
&\le C \|G_{\tau,h}f - \Pi_h(Gf)\|_{L^2(\Omega)}^2 \\
&\le C \|G_{\tau,h}f - Gf\|_{L^2(\Omega)}^2 + C\|\Pi_h(Gf) - Gf\|_{L^2(\Omega)}^2 \\
&\le C \|G_{\tau,h}f - Gf\|_{L^2(\Omega)}^2 + Ch^4\|Gf\|_{H^2(\Omega)}^2.
\end{align*}
Combining these parts gives the bound
$\| Gf - G_{\tau,h} f \|_n^2 \le C \| G_{\tau,h} f - Gf \|_{L^2(\Omega)}^2 + C h^4 \| Gf \|_{H^2(\Omega)}^2.$
To estimate the first term on the right, we apply the standard error estimate for the numerical scheme \eqref{5.2} from \cite[Theorem 4.1]{35}. As established in Lemma \ref{L5.1}, our solution possesses $H^2$ spatial regularity, so we have
\begin{align*}
\| Gf - G_{\tau,h} f \|_{L^2(\Omega)} \le C \bigg(
& h^2 \left( \| u \|_{L^\infty(0,T;H^2(\Omega))} + \| u_t \|_{L^2(0,T;H^2(\Omega))} \right) \\
+& \tau^2 \left( \| u_{ttt} \|_{L^2(0,T;L^2(\Omega))} + \| u_{tttt} \|_{L^2(0,T;L^2(\Omega))} \right)
\bigg).
\end{align*}
Invoking the regularity estimates from Lemma \ref{L5.1}, all norms on the right-hand side are bounded by $C\|f\|_{L^2(\Omega)}$. Substituting this back concludes the proof.
\end{proof}

\section{Convergence Analysis for the Continuous Solution}\label{2}
This 
section
 presents the  convergence analysis for the regularized solutions in the continuous setting under noise models (Y1) and (Y2).

\subsection{Analysis under Noise Model (Y1)}

The proof of the main theorem relies on the characterization of the solution space $V_n$ and the spectral bounds of the forward operator, as established in \cite[Lemmas 2.1 and 2.2]{1}.

\begin{theorem}\label{Ta3.3}
Suppose Lemma~\ref{T3.1} holds, and let $f_n \in L^2(\Omega)$ be the solution of \eqref{2.2}. There exist $\lambda_0 > 0$ and $C > 0$ such that for all $\alpha \leq \lambda_0$,
\begin{align}\label{3.4a}
    \mathbb{E} \left[ \|G f_n - G f^*\|_n^2 \right] &\leq C \alpha \|f^*\|_{L^2(\Omega)}^2 + C \sigma^2 (n \alpha^{d/4})^{-1}, \\\label{3.5a}
    \mathbb{E} \left[ \|f_n - f^*\|_{L^2(\Omega)}^2 \right] &\leq C \|f^*\|_{L^2(\Omega)}^2 + C \sigma^2 (n \alpha^{1+d/4})^{-1}.
\end{align}
\end{theorem}
\begin{proof}
The minimizer $f_n$ of \eqref{2.2} satisfies the variational equation \begin{equation}\label{3.3ab}
\alpha (f_n, v)_{L^2(\Omega)} + (G f_n, G v)_n = (m, G v)_n, ~~for ~~all~~ v~ \in L^2(\Omega).\end{equation}
Define the norm $\|v\|_{\alpha}^2 := \alpha \|v\|_{L^2(\Omega)}^2 + \|G v\|_n^2$. Setting $v=f_n-f^*$ and using \eqref{2.3} gives
\begin{equation}\label{3.7a}
\begin{split}
\|f_n - f^*\|_{\alpha}^2 &= (e, G(f_n-f^*))_n - \alpha(f^*, f_n-f^*)_{L^2(\Omega)} \\
&\le \alpha \|f^*\|_{L^2(\Omega)}^2 + \sup_{v \in L^2(\Omega)} \frac{(e, G v)^2_n}{\|v\|_{\alpha}^2}.
\end{split}
\end{equation}

Since the solution to the optimization problem 
$\min_{f \in L^2(\Omega), (Gf)(x_i) = m_i} \|f\|_{L^2(\Omega)}^2 $lies in an $n$-dimensional subspace $V_n$ \cite[Lemma 2.1]{1}, we can restrict the supremum as follows:
\begin{align*}
\sup_{v \in L^2(\Omega)} \frac{(e, Gv)_n^2}{||v||_{\alpha}^2} &= \sup_{v \in L^2(\Omega)} \frac{(e, Gv)_n^2}{\alpha\|v\|_{L^2(\Omega)}^2 + \|Gv\|_n^2} \\
&\le \sup_{v \in L^2(\Omega)} \frac{(e, Gv)_n^2}{\alpha \min_{u \in L^2(\Omega), Gu(x_i)=Gv(x_i)} \|u\|_{L^2(\Omega)}^2 + \|Gv\|_n^2} \\
&= \sup_{v \in L^2(\Omega)} \frac{(e, Gv)_n^2}{\alpha \min_{u \in V_n, Gu(x_i)=Gv(x_i)} \|u\|_{L^2(\Omega)}^2 + \|Gv\|_n^2} \\
&= \sup_{v \in V_n} \frac{(e, Gv)_n^2}{||v||_{\alpha}^2}.
\end{align*}

According to \cite[Lemma 2.2]{1}, we can let $\rho_1 \le \rho_2 \le \cdots \le \rho_n$ be the eigenvalues of the problem
$
(\phi, v)_{L^2(\Omega)} = \rho(G\phi, Gv)_n,  \forall v \in V_n$,
with the corresponding eigenfunctions $\{\phi_k\}_{k=1}^n$, which is an orthonormal basis of $V_n$ under the inner product $(G\cdot, G\cdot)_n$. Thus $(G\phi_k, G\phi_l)_n = \delta_{kl}$ and consequently, $(\phi_k, \phi_l)_{L^2(\Omega)} = \rho_k \delta_{kl}$, $l=1,2,\dots,n$.

Now for any $v \in V_n$, we have the expansion $v(x) = \sum_{k=1}^n v_k \phi_k(x)$, where $v_k = (Gv, G\phi_k)_n$ for $k=1,2,\dots,n$. Thus $||v||_{\alpha}^2 = \sum_{k=1}^n (\alpha \rho_k + 1)v_k^2$. By the Cauchy-Schwarz inequality we can readily get
\begin{align*}
(e, Gv)_n^2 &= \frac{1}{n^2} \left( \sum_{i=1}^n e_i \left( \sum_{k=1}^n v_k G\phi_k(x_i) \right) \right)^2 = \frac{1}{n^2} \left( \sum_{k=1}^n v_k \left( \sum_{i=1}^n e_i G\phi_k(x_i) \right) \right)^2 \\
&\le \frac{1}{n^2} \sum_{k=1}^n (1+\alpha\rho_k)v_k^2 \cdot \sum_{k=1}^n (1+\alpha\rho_k)^{-1} \left( \sum_{i=1}^n e_i (G\phi_k)(x_i) \right)^2.
\end{align*}

Taking the expectation and using the independence of $\{e_i\}$ with $\mathbb{E}[e_i^2] \leq \sigma^2$, we have 
\begin{align*}
\mathbb{E}\left[ \sup_{v \in V_n} \frac{(e, G v)_n^2}{\|v\|_{\alpha}^2} \right] &\le \frac{1}{n^2} \sum_{k=1}^n \frac{1}{\alpha \rho_k + 1} \mathbb{E}\left[\left( \sum_{i=1}^n e_i (G\phi_k)(x_i) \right)^2\right] \\
&= \frac{\sigma^2}{n^2} \sum_{k=1}^n \frac{1}{\alpha \rho_k + 1} \sum_{i=1}^n (G\phi_k)(x_i)^2 = \frac{\sigma^2}{n} \sum_{k=1}^n \frac{\|G\phi_k\|_n^2}{\alpha \rho_k + 1}.
\end{align*}
Since $\|G\phi_k\|_n=1$ and using the eigenvalue bound $\rho_k \ge Ck^{4/d}$ from \cite[Lemma 2.2]{1}, we can approximate the sum by an integral: $\sum_{k=1}^n \frac{1}{\alpha \rho_k + 1} \le \int_0^\infty (1+C\alpha t^{4/d})^{-1}dt \le C' \alpha^{-d/4}$. Thus,
\[
\mathbb{E} \left[ \sup_{v \in V_n} \frac{(e, G v)_n^2}{\|v\|_{\alpha}^2} \right] \leq C'' \sigma^2 (n \alpha^{d/4})^{-1}.
\]
Combining this with \eqref{3.7a} gives $\mathbb{E}[\|f_n - f^*\|_{\alpha}^2] \leq \alpha \|f^*\|_{L^2(\Omega)}^2 + C'' \sigma^2 (n \alpha^{d/4})^{-1}$. The bounds \eqref{3.4a} and \eqref{3.5a} follow from the definition of the $\|\cdot\|_{\alpha}$ norm.
\end{proof}


Having established the $L^2$-norm convergence, we now derive the corresponding convergence rate in the weaker $H^{-1}(\Omega)$ 
norm.
\begin{theorem}\label{T3.3abc}
Suppose Theorem \ref{Ta3.3} holds, and $\alpha \geq n^{-4/d}$. Then
\[
\mathbb{E}\left[\|f_n - f^*\|_{H^{-1}(\Omega)}^2\right] \leq C\alpha^{1/2}\|f^*\|_{L^2(\Omega)}^2 + C\sigma^2(n\alpha^{1/2 + d/4})^{-1}.
\]
\end{theorem}
\begin{proof}
By \eqref{2.7ab}, for any $v \in L^2(\Omega)$, we have the expansion $v = \sum_{k=1}^{\infty} v_k \psi_k$ with $v_k = (v, \psi_k)_{L^2(\Omega)}$. We can directly check that $\|v\|_{L^2(\Omega)}^2 = \sum_{k=1}^{\infty} v_k^2$ and $\|Gv\|_{L^2(\Omega)}^2 = \sum_{k=1}^{\infty} \eta_k^{-1} v_k^2$. Then for any $g \in L^2(\Omega)$, $g = \sum_{k=1}^{\infty} g_k \phi_k$, with $g_k = (g, \psi_k)_{L^2(\Omega)}$, we can obtain by the Cauchy-Schwarz inequality that
\begin{align*}
\|v\|_{H^{-1}(\Omega)} &= \sup_{0\neq g \in H_0^1(\Omega)} \frac{|(g,v)_{L^2(\Omega)}|}{\|g\|_{H_0^1(\Omega)}} = \sup_{0\neq g \in H_0^1(\Omega)} \frac{|\sum_{k=1}^\infty g_k v_k|}{\|g\|_{H_0^1(\Omega)}} \\
&\le \left( \sum_{k=1}^\infty \eta_k^{-1/2} v_k^2 \right)^{1/2} \le \|Gv\|_{L^2(\Omega)}^{1/2} \|v\|_{L^2(\Omega)}^{1/2}.
\end{align*}

Taking $v = f^* - f_n$ in the above inequality, we obtain 
\begin{equation}\label{3.12a_full}
\|f^* - f_n\|_{H^{-1}(\Omega)}^2 \leq C \|G(f^* - f_n)\|_{L^2(\Omega)} \|f^* - f_n\|_{L^2(\Omega)}.
\end{equation}
The first term on the right-hand side can be bounded by relating the continuous and discrete norms. Using the norm relationship from \cite[Theorems 3.3 and 3.4]{39}, which for any $u \in H^2(\Omega)$ is given by
$\|u\|_{L^2(\Omega)}^2 \leq C \left( \|u\|_n^2 + n^{-4/d} \|u\|_{H^2(\Omega)}^2 \right)$,
together with the condition $\alpha \geq n^{-4/d}$ and the boundedness of the operator $G: L^2(\Omega)\to H^2(\Omega)$, we obtain the key intermediate bound:
\begin{equation}\label{3.15ab}
\|G(f^* - f_n)\|_{L^2(\Omega)}^2 \le C \left( \|G(f^* - f_n)\|_n^2 + \alpha \|f^* - f_n\|_{L^2(\Omega)}^2 \right) = C\|f^*-f_n\|_{\alpha}^2.
\end{equation}
Applying Young's inequality to \eqref{3.12a_full} yields:
\begin{equation}\label{3.15a}
\|f^* - f_n\|_{H^{-1}(\Omega)}^2 \leq C \left( \alpha^{-1/2} \|G(f^* - f_n)\|_{L^2(\Omega)}^2 + \alpha^{1/2} \|f^* - f_n\|_{L^2(\Omega)}^2 \right).
\end{equation}
Combining these gives $\|f^*-f_n\|_{H^{-1}(\Omega)}^2 \le C\alpha^{-1/2}\|f^*-f_n\|_{\alpha}^2$. Taking the expectation and applying Theorem \ref{Ta3.3}, we conclude:
\begin{align*}
\mathbb{E}\left[\|f_n - f^*\|_{H^{-1}(\Omega)}^2\right] &\le C\alpha^{-1/2} \mathbb{E}\left[\|f_n-f^*\|_{\alpha}^2\right] \\
&\le C\alpha^{-1/2} \left( \alpha \|f^*\|_{L^2(\Omega)}^2 + C \sigma^2 (n\alpha^{d/4})^{-1} \right) \\
&= C \alpha^{1/2} \|f^*\|_{L^2(\Omega)}^2 + C \sigma^2 (n \alpha^{1/2 + d/4})^{-1}.
\end{align*}
\end{proof}

\subsection{
Analysis under Noise Model (Y2)}

This subsection is devoted to the convergence analysis under the sub-Gaussian noise model (Y2
), which assumes that the noise variables satisfy
\[
\mathbb{E} \left[ \exp\left( \lambda (e_i - \mathbb{E}[e_i]) \right) \right] \leq \exp\left( \frac{1}{2} \sigma^2 \lambda^2 \right), \quad \forall \lambda \in \mathbb{R}.
\]
We analyze the  convergence of the error terms $\|Gf^* - Gf_n\|_n$ and $\|f^* - f_n\|_{H^{-1}(\Omega)}$.

\subsubsection{Preliminaries on Sub-Gaussian Theory}

We first recall essential concepts from sub-Gaussian theory that form the foundation of our analysis (see \cite{12,13,14}). A random variable $Z$ is called sub-Gaussian if its tail probabilities exhibit exponential decay, i.e., $\mathbb{P}(|Z - \mathbb{E}[Z]| \geq z) \leq 2 \exp( -z^2/(2\sigma^2))$. An equivalent characterization is provided by the Orlicz norm. For $\psi_2(t) = e^{t^2} - 1$, the Orlicz norm $\|Z\|_{\psi_2}$ is defined
 as
\[
\|
Z\|_{\psi_2} = \inf \left\{ C > 0 : \mathbb{E} \left[ \psi_2 \left( \frac{|Z|}{C} \right) \right] \leq 1
 \right\},
\]
which yields the tail 
bound:
\begin{equation}\label{4.3a}
P(|Z| \geq z) \leq 2 \exp\left( -z^2/\|Z\|_{\psi_2}^2 \right) \quad \forall z > 0.
\end{equation}

A random process $\{Z_t : t \in \mathbb{T}\}$ indexed by a semi-distance space $(\mathbb{T}, \hat{\mathbf{d}})$ is sub-Gaussian if $P(|Z_s - Z_t| > z) \leq 2 \exp(-z^2/(2 \hat{\mathbf{d}}(s,t)^2))$. The complexity of the index set $\mathbb{T}$ is measured by its covering number, $N(\varepsilon, \mathbb{T},\hat{\mathbf{d}})$, defined as the minimum number of $\varepsilon$-balls in the semi-distance $\hat{\mathbf{d}}$ required to cover $\mathbb{T}$.

\subsubsection{Convergence Results under Sub-Gaussian Noise}
Equipped with these probabilistic tools, we now establish convergence rates in probability for
 the continuous regularized solution.
\begin{theorem}\label{T4.1}
Suppose Lemma~\ref{T3.1} holds, let $\rho_0 = \|f^*\|_{L^2(\Omega)} + \sigma n^{-1/2}$, and let $f_n$ be the solution of \eqref{2.2}. Under the condition $\alpha^{1/2 + d/8} = O(\sigma n^{-1/2} \rho_0^{-1})$, there exists a constant $C > 0$ such that for $z \ge 1$,
\begin{align*}
\mathbb{P}\left(\|Gf_n - Gf^*\|_n \geq \alpha^{1/2} \rho_0 z\right) &\leq 2e^{-Cz^2}, \\
\mathbb{P}\left(\|f_n\|_{L^2(\Omega)} \geq \rho_0 z\right) &\leq 2e^{-Cz^2}.
\end{align*}
\end{theorem}
\begin{proof}
The proof employs the peeling technique. The optimality condition of \eqref{2.2} yields:
\begin{equation}\label{4.6a}
\|G f_n - G f^*\|_n^2 + \alpha \|f_n\|_{L^2(\Omega)}^2 \leq 2 (e, G f_n - G f^*)_n + \alpha \|f^*\|_{L^2(\Omega)}^2.
\end{equation}
Let \(\delta\), \(\rho > 0\) be two constants to be determined later, and let \(i, j \geq 1\)
\begin{equation*}
A_0 = [0, \delta), \quad A_i = [2^{i-1} \delta, 2^i \delta), \quad B_0 = [0, \rho), \quad B_j = [2^{j-1} \rho, 2^j \rho).
\end{equation*}
For $i, j \geq 0$, consider the sets $F_{ij} = \{ v \in L^2(\Omega) : \| G v -G f^*\|_n \in A_i, \| v \|_{L^2(\Omega)} \in B_j \}$. This gives:
\begin{equation}\label{4.8a}
\mathbb{P}(\| G f_n - G f^* \|_n > \delta) \leq \sum_{i=1}^{\infty} \sum_{j=0}^{\infty} \mathbb{P}(f_n  \in F_{ij}).
\end{equation}

By \cite[Lemma 4.6]{12}, the process $\{ (e, G v)_n : v \in L^2(\Omega) \}$ is sub-Gaussian with respect to the semi-distance $\hat{\mathbf{d}}(f, v) = \sigma n^{-1/2} \| G f - G v \|_n$. With this semi-distance, it is easy to see that $\text{diam}(F_{ij}) \le 2\sigma n^{-1/2} \cdot 2^i \delta$. Using Dudley's entropy integral theorem for sub-Gaussian processes \cite{14} then implies:
\begin{align*}
\sup_{f \in F_{ij}} \| (e, Gf - Gf^*)_n \|_{\Psi_2} &\le K \int_0^{\sigma n^{-1/2} \cdot 2^{i+1} \delta} \sqrt{\log N\left(\frac{\epsilon}{2}, F_{ij}, \hat{\mathbf{d}}\right)} \, d\varepsilon \\
&= K \int_0^{\sigma n^{-1/2} \cdot 2^{i+1} \delta} \sqrt{\log N\left(\frac{\epsilon}{2\sigma n^{-1/2}}, F_{ij}, \|G \cdot \|_n\right)} \, d\epsilon.
\end{align*}
By Lemma \ref{L5.1}, the range of the forward operator $G$ is in $H^2(\Omega)$, which allows us to apply the  Kolmogorov's estimate for the metric entropy of Sobolev spaces  from \cite{40} to obtain the covering entropy bound:
\begin{equation*}
\begin{split}
\log N \left( \frac{\epsilon}{2 \sigma n^{-1/2}}, F_{ij}, \| G \cdot \|_n \right) &\le \log N \left( \frac{\epsilon}{2 \sigma n^{-1/2}}, G(F_{ij}), \| \cdot \|_{L^\infty(\Omega)} \right) \\&\leq C \left( \frac{2 \sigma n^{-1/2} \cdot 2^j \rho}{\epsilon} \right)^{d/2},
\end{split}
\end{equation*}
where the last inequality uses that $G(F_{ij})$ is contained in a ball in $H^2(\Omega)$ of radius $C(2^j \rho)$, since $G : L^2(\Omega) \to H^2(\Omega)$ is bounded. Consequently:
\begin{equation*}
\| \sup_{f - f^* \in F_{ij}} |(e, G f - G f^*)_n| \|_{\psi_2} \leq C \sigma n^{-1/2} (2^j \rho)^{d/4} (2^i \delta)^{1-d/4}.
\end{equation*}
For $i, j \geq 1$, combining \eqref{4.6a} with the tail bound \eqref{4.3a} yields:
\begin{align*}
\mathbb{P}(f_n  \in F_{ij})&\leq \mathbb{P} \left(   2^{2(i-1)} \delta^2 + \alpha 2^{2(j-1)} \rho^2\leq 2 \sup_{f - f^* \in F_{ij}} |(e, G f - G f^*)_n| + \alpha \rho_0^2\right) \\ &\leq \mathbb{P} \left( 2 \sup_{f - f^* \in F_{ij}} |(e, G f - G f^*)_n| \geq 2^{2(i-1)} \delta^2 + \alpha 2^{2(j-1)} \rho^2 - \alpha \rho_0^2 \right) \\
&\leq 2\exp \left[ -\frac{1}{C\sigma^2 n^{-1}} \left( \frac{2^{2(i-1)} \delta^2 + \alpha 2^{2(j-1)} \rho^2 - \alpha \rho_0^2}{(2^i \delta)^{1-d/4} (2^j \rho)^{d/4}} \right)^2 \right].
\end{align*}
Set $z \geq 1$, $\delta^2 = \alpha \rho_0^2 (1 + z)^2$, $\rho = \rho_0$. Under the assumed condition on $\alpha$, direct computation gives for $i,j \geq 1$:
\begin{equation}\label{4.10a}
\mathbb{P}(f_n  \in F_{ij}) \leq 2\exp \left[ -C \left( \frac{2^{2(i-1)} z(1 + z) + 2^{2(j-1)}}{ (2^i (1 + z))^{1 - d/4} (2^j )^{d/4} } \right)^2 \right].
\end{equation}
By Young's inequality, the denominator is bounded by $C ((1 + z)2^i + 2^j)$, which simplifies the bound to:
\[
\mathbb{P}(f_n\in F_{ij}) \leq 2\exp \left[ -C(2^{2i} z^2 + 2^{2j}) \right].
\]
A similar argument for $i \geq 1, ~j = 0$ gives $\mathbb{P}(f_n \in F_{i0}) \leq 2 \exp [ -C (2^i z^2) ]$. Summing the probabilities and using $\sum_{k=1}^\infty \exp(-C2^{2k}) \leq\exp(-C)< 1$ and $\sum_{k=1}^\infty \exp(-C2^{2k}z^2) < \exp(-Cz^2)$:
\[
\sum_{i=1}^\infty \sum_{j=0}^\infty \mathbb{P}(f_n \in F_{ij}) \le 2 \sum_{i=1}^\infty \sum_{j=1}^\infty \exp(-C(2^{2i}z^2 + 2^{2j})) + 2 \sum_{i=1}^\infty \exp(-C(2^{2i}z^2)).
\]

From \eqref{4.8a}, this implies $\mathbb{P}(\| G f_n - G f^* \|_n > \alpha^{1/2} \rho_0 (1 + z)) \leq 4 \exp(-C z^2)$ for $z \geq 1$. By \cite[Lemma 4.7]{12}, this establishes:
\begin{equation}\label{4.5a}
\|\| G f_n - G f^* \|_n \|_{\psi_2} \leq C \alpha^{1/2} \rho_0.
\end{equation}
Similar arguments (summing over $i \geq 0$, $j \geq 1$) yield:
\begin{equation}\label{4.5b}
\|\| f_n \|_{\psi_2} \leq C \rho_0.
\end{equation}
Applying the Orlicz norm property \eqref{4.3a} to \eqref{4.5a} and \eqref{4.5b} gives the desired probability bounds.
\end{proof}

Building upon Theorem~\ref{T4.1}, we derive the corresponding convergence rate in the $H^{-1}$ 
topology.
\begin{theorem}\label{T4.2}
Suppose Theorem \ref{T4.1} holds, and assume $\alpha \geq n^{-4/d}$. Then for $z \ge 1$,
\[
\mathbb{P}\left(\|f_n - f^*\|_{H^{-1}(\Omega)} \geq \alpha^{1/4}\rho_0 z\right) \leq 2e^{-Cz^2}.
\]
\end{theorem}
\begin{proof}
Starting from inequalities \eqref{3.15ab} and \eqref{3.15a},  we have:
\begin{equation}\label{step1}
\|f^* - f_n\|_{H^{-1}(\Omega)}^2 \leq C \left( \alpha^{1/2} \|f^* - f_n\|_{L^2(\Omega)}^2 + \alpha^{-1/2} \|G(f^* - f_n)\|_n^2 \right).
\end{equation}
Taking the Orlicz norm and applying the triangle inequality yields:
\begin{equation}\label{step2}
\| \|f^* - f_n\|_{H^{-1}(\Omega)} \|_{\psi_2} \leq C \alpha^{1/4} \| \|f^* - f_n\|_{L^2(\Omega)} \|_{\psi_2} + C \alpha^{-1/4} \| \|G(f^* - f_n)\|_n \|_{\psi_2}.
\end{equation}
From inequalities \eqref{4.5a} and \eqref{4.5b}, we have $\| \|G f_n - G f^*\|_n \|_{\psi_2} \leq C \alpha^{1/2} \rho_0$ and $\| \|f_n \| \|_{\psi_2} \leq C \rho_0$. This also implies $\| \|f^*-f_n\| \|_{\psi_2} \le \|\|f^*\|\|_{\psi_2} + \|\|f_n\|\|_{\psi_2} \le C\rho_0$. Substituting these into \eqref{step2} gives:
\begin{equation}\label{step3}
\| \|f^* - f_n\|_{H^{-1}(\Omega)} \|_{\psi_2} \leq C \alpha^{1/4} \rho_0 + C \alpha^{-1/4} \cdot (\alpha^{1/2} \rho_0) = C \rho_0 \alpha^{1/4}.
\end{equation}
The result follows from the Orlicz norm property \eqref{4.3a} and the bound in \eqref{step3}.
\end{proof}

\section{Convergence Analysis for the Discrete Solution}\label{3}
This section establishes error estimates for the fully discrete regularized solution under both noise models (Y1) and (Y2).

\subsection{Convergence under Noise Model (Y1)}

We now present the main convergence result for the fully discrete estimator $f_h$ defined in \eqref{5.6}.

\begin{theorem}\label{T5.1}
Assume Lemma~\ref{T3.1} holds, let $f_n$ solve \eqref{2.2} and $f_h \in V_h$ solve \eqref{5.6}. Assume $\tau = O(h)$. There exist $\lambda_0 > 0, C>0$ such that for $\alpha \le \lambda_0$ and $h^4 \le C\alpha$,
\begin{align}
    \mathbb{E}\left[ \| G f^* - G_{\tau, h} f_h \|_n^2 \right] &\leq C \alpha \| f^* \|_{L^2(\Omega)}^2 + C \frac{\sigma^2}{n \alpha^{d/4}}, \label{5.22} \\
    \mathbb{E}\left[ \| f^* - f_h \|_{L^2(\Omega)}^2 \right] &\leq C \| f^* \|_{L^2(\Omega)}^2 + C \frac{\sigma^2}{n \alpha^{1+d/4}}. \label{5.23}
\end{align}
\end{theorem}
\begin{proof}
Let's define the bilinear form $a_h(f, v) = \alpha (f, v)_{L^2(\Omega)} + (G_{\tau,h} f, G_{\tau,h} v)_n$ and the associated energy norm $\|f\|_{a_h}^2 = a_h(f, f)$ for any $f, v \in L^2(\Omega)$. Taking $v=w_h$ in \eqref{3.3ab} and $v_h=w_h$ in \eqref{2.11a}, 
for any $v_h, w_h \in V_h$, we can derive the error equation:
\[
a_h(f_n - v_h, w_h) = a_h(f_n-v_h, w_h) + F(w_h),
\]
where the consistency error functional $F(w_h)$ is given by
\[
F(w_h) = ((G - G_{\tau,h}) f_n, G_{\tau,h} w_h)_n + (G f^* - G f_n, (G_{\tau,h} - G) w_h)_n + (e, (G_{\tau,h} - G) w_h)_n.
\]
Using the triangle inequality, we have:
\begin{equation}\label{3.19a_full}
\|f_n - f_h\|_{a_h} \leq C \inf_{v_h \in V_h} \|f_n - v_h\|_{a_h} + C \sup_{w_h \in V_h} \frac{|F(w_h)|}{\|w_h\|_{a_h}}.
\end{equation}
We now bound the expectation of each error component. 
The two deterministic terms in $F(w_h)$ are bounded using Lemma \ref{L2.3a}:
\begin{align}
\sup_{w_h\in V_h} \frac{|((G - G_{\tau,h}) f_n, G_{\tau,h} w_h)_n|}{\|w_h\|_{a_h}} &\leq \|G f_n - G_{\tau,h} f_n\|_n \leq C h^2 \|f_n\|_{L^2(\Omega)}, \label{3.20a_full} \\
\sup_{w_h\in V_h} \frac{|(G f^* - G f_n, (G_{\tau,h} - G) w_h)_n|}{\|w_h\|_{a_h}} &\leq C \|G f^* - G f_n\|_n \frac{h^2}{\alpha^{1/2}}. \label{3.21a_full}
\end{align}
For the  random term, let $\{\xi_j\}_{j=1}^{N_h}$ be an $L^2$-orthonormal basis for the finite element space $V_h$, where $N_h = \dim(V_h)$, such that $(\xi_i, \xi_j)_{L^2(\Omega)} = \delta_{ij}$. For any $w_h = \sum_j (w_h, \xi_j)\xi_j$, the Cauchy-Schwarz inequality gives
\begin{align*}
(e, (G - G_{\tau,h}) w_h)_n^2 &\le \frac{1}{n^2} \sum_{j=1}^{N_h} (w_h, \xi_j)^2 \sum_{j=1}^{N_h} \left( \sum_{i=1}^n e_i (G - G_{\tau,h}) \xi_j(x_i) \right)^2 \\
&= \frac{1}{n^2} \|w_h\|_{L^2(\Omega)}^2 \sum_{j=1}^{N_h} \left( \sum_{i=1}^n e_i (G - G_{\tau,h}) \xi_j(x_i) \right)^2,
\end{align*}
Dividing by $\|w_h\|_{a_h}^2 \ge \alpha \|w_h\|_{L^2(\Omega)}^2$ and using $\mathbb{E}[e_i^2] \leq \sigma^2$ and Lemma \ref{L2.3a}, we get
\begin{align*}
\mathbb{E}\left( \sup_{w_h \in V_h} \frac{|(e, G w_h - G_{\tau,h} w_h)_n|^2}{\|w_h\|_{_{a_h}}^2} \right) &\le \frac{1}{\alpha n^2} \sum_{j=1}^{N_h} \mathbb{E}\left( \sum_{i=1}^n e_i (G - G_{\tau,h})\xi_j(x_i) \right)^2 \\
&\le \frac{1}{\alpha n} \sum_{j=1}^{N_h} \sigma^2 \|(G - G_{\tau,h})\xi_j\|_n^2 \le C \frac{\sigma^2}{\alpha n} N_h h^4.
\end{align*}
Since the $L^2$-orthogonal projection $P_h : L^2(\Omega) \to V_h$ satisfies $G_{\tau,h}f = G_{\tau,h}(P_h f)$ by \eqref{5.2}, we obtain the approximation error estimate
$
\|f_n - P_h f_n\|_{a_h}^2 \le C(\alpha + h^4) \|f_n\|_{L^2(\Omega)}^2.
$
Combining this estimate with the bounds on the consistency terms, we arrive at:
\begin{align*}
\begin{split}
\mathbb{E}\left[ \|f_n - f_h\|_{a_h}^2 \right] \leq& C\alpha\mathbb{E}[\|f_n\|_{L^2(\Omega)}^2] + C\frac{h^4}{\alpha}\mathbb{E}[\|Gf^*-Gf_n\|_n^2] \\&+ C \frac{\sigma^2 N_h h^4}{\alpha n}.
\end{split}
\end{align*}
Using the bounds on the continuous error from Theorem \ref{Ta3.3} and the property $N_h \le Ch^{-d}$ for quasi-uniform meshes,  we find
\begin{equation}
\begin{split}\label{3.13ab}
\mathbb{E}\left[ \| G f^* - G_{\tau, h} f_h \|_n^2 \right] &\leq C (\alpha + h^4) \| f^* \|_{L^2(\Omega)}^2 + C \left( 1 + \frac{h^4}{\alpha} \right) \frac{\sigma^2}{n \alpha^{d/4}}, \\
\mathbb{E}\left[ \| f^* - f_h \|_{L^2(\Omega)}^2 \right] &\leq C \frac{\alpha + h^4}{\alpha} \| f^* \|_{L^2(\Omega)}^2 + C \left( 1 + \frac{h^4}{\alpha} \right) \frac{\sigma^2}{n \alpha^{1+d/4}}.
\end{split}
\end{equation}
The final result follows from the triangle inequality and the assumption $h^4 \leq C\alpha$.
\end{proof}

Finally, we translate the discrete $L^2$-norm convergence into a convergence rate in the weaker $H^{-1}$ 
topology.

\begin{theorem}\label{T3.7a}
Suppose Theorem~ \ref{T5.1} holds, and $\alpha \geq n^{-4/d}$. Then
\begin{align}\label{5.24}
    \mathbb{E}\left[ \| f^* - f_h \|_{H^{-1}(\Omega)}^2 \right] \leq C \alpha^{1/2}   \| f^* \|_{L^2(\Omega)}^2 + C  \frac{\sigma^2}{n \alpha^{d/4+1/2}}.
\end{align}
\end{theorem}
\begin{proof}
Following the proof of \eqref{3.15ab} and \eqref{3.15a}, we have
\begin{equation}\label{step1a}
\| f^* - f_h \|_{H^{-1}(\Omega)}^2 \leq C \alpha^{1/2} \|f^* - f_h\|_{L^2(\Omega)}^2 + C \alpha^{-1/2} \|G(f^* - f_h)\|_{n}^2.
\end{equation}
By the triangle inequality, $\|G(f^*-f_h)\|_{n} \le \|Gf^*-G_{\tau,h}f_h\|_{n} + \|G_{\tau,h}f_h - Gf_h\|_{n}$. By Lemma \ref{L2.3a} and $\tau = O(h)$,
\begin{equation}\label{step2a}
\| G_{\tau,h} f_h - G f_h \|_n^2 \leq C (h^2 + \tau^2)^2 \| f_h \|_{L^2(\Omega)}^2 = C h^4 \| f_h \|_{L^2(\Omega)}^2.
\end{equation}
Taking the expectation of \eqref{step1a} and substituting the bounds gives:
\begin{equation}\label{step3a}
\begin{split}
\mathbb{E}\left[ \| f^* - f_h \|_{H^{-1}(\Omega)}^2 \right] \leq& C \alpha^{1/2} \mathbb{E}\left[ \|f^* - f_h\|_{L^2(\Omega)}^2 \right] + C \alpha^{-1/2} \mathbb{E}\left[ \|G f^* - G_{\tau,h} f_h\|_n^2 \right] \\&+ C \alpha^{-1/2} h^4 \mathbb{E}\left[ \| f_h \|_{L^2(\Omega)}^2 \right].
\end{split}
\end{equation}
We bound $\mathbb{E}[\|f_h\|_{L^2(\Omega)}^2] \le 2\mathbb{E}[\|f^*-f_h\|_{L^2(\Omega)}^2] + 2\|f^*\|_{L^2(\Omega)}^2$. Substituting the bounds from Theorem \ref{T5.1} for each term and collecting terms leads to the desired estimate.
\end{proof}

\subsection{Convergence under Noise Model (Y2)}

This subsection analyzes the convergence of the discrete solution \eqref{5.6} under sub-Gaussian noise. The following lemma establishes bounds on the Orlicz norm of the random component in the consistency error, which is essential for the discrete analysis.

\begin{lemma}\label{L4.6}
Let $B_h := \{ w_h \in V_h : \|w_h\|_{a_h} \leq 1 \}$, where $V_h$ is the finite element space with dimension $N_h = \dim(V_h)$ satisfying $N_h \le C h^{-d}$ for a quasi-uniform mesh. The energy norm is defined by $\|w_h\|_{a_h}^2 = \alpha\|w_h\|_{L^2(\Omega)}^2+\|G_{\tau,h}w_h\|_n^2$. Assume that $h^4 \leq C \alpha$. Then the Orlicz norm of the supremum of the discretization error process is bounded by:
\[
\left\| \sup_{w_h \in B_h} \left| (e, G w_h - G_{\tau,h} w_h)_n \right| \right\|_{\psi_2} \leq C \sigma n^{-1/2} \alpha^{-d/8}.
\]
\end{lemma}
\begin{proof}
 Let $\widehat{E}_n(w_h) := (e, (G - G_{\tau,h}) w_h)_n$. This process is known to be sub-Gaussian (see \cite[Lemma 4.6]{12}), and its semi-distance can be bounded using the discretization error estimate from Lemma~\ref{L2.3a} and \cite[Lemma 4.6]{12}:
$\hat{\mathbf{d}}(v_h, w_h) = \sigma n^{-1/2} \|(G v_h - G_{\tau,h} v_h) - (G w_h - G_{\tau,h} w_h)\|_n \le C \sigma n^{-1/2} h^2 \|v_h - w_h\|_{L^2(\Omega)}$.
For any $w_h \in B_h$, the definition of the $\|\cdot\|_{a_h}$ norm implies that $\|w_h\|_{L^2(\Omega)} \le \alpha^{-1/2}$. The diameter of $B_h$ is therefore bounded by $D_h := \sup_{v_h, w_h \in B_h} \hat{\mathbf{d}}(v_h, w_h) \le C \sigma n^{-1/2} h^2 \alpha^{-1/2}$.

Dudley's entropy integral for sub-Gaussian processes  \cite{14} gives \begin{equation}\label{4.10ab}
\|\sup_{w_h \in B_h} |\widehat{E}_n(w_h)|\|_{\psi_2} \le K \int_0^{D_h} \sqrt{\log N(\epsilon/2, B_h, \hat{\mathbf{d}})} \, d\epsilon.\end{equation}
The covering number can be bounded using a standard result for finite-dimensional spaces \cite[Corollary 2.6]{13}:
\[
\log N \left( \epsilon/2, B_h, \hat{\mathbf{d}} \right) \leq \log N \left( \frac{\epsilon}{C \sigma n^{-1/2} h^2}, B_h, \| \cdot \|_{L^2(\Omega)} \right) \leq N_h \left( 1 + \frac{C' D_h}{\epsilon} \right).
\]
By applying the Dudley entropy integral theorem, we obtain the bound $C \sqrt{N_h} D_h$ in \eqref{4.10ab}. Substituting the expression for the diameter $D_h$, we find:
\begin{align*}
\left\| \sup_{w_h \in B_h} |\widehat{E}_n(w_h)| \right\|_{\psi_2} &\leq C \sqrt{N_h} \left( \sigma n^{-1/2} h^2 \alpha^{-1/2} \right) \\
&= C \sigma n^{-1/2} \alpha^{-1/2} \sqrt{N_h h^4}.
\end{align*}
Using the property $N_h \le C h^{-d}$ and the assumption $h^4 \le C\alpha$ (which implies $O(h^4) \le C\alpha$, so $h \le C\alpha^{1/4}$), we have:
\[
N_h h^4 \le  O(h^{4-d}) \le C (\alpha^{1/4})^{4-d} = C \alpha^{1-d/4}.
\]
Substituting this derived bound back into our expression gives:
\begin{align*}
\left\| \sup_{w_h \in B_h} |\widehat{E}_n(w_h)| \right\|_{\psi_2} &\leq C \sigma n^{-1/2} \alpha^{-1/2} \sqrt{C \alpha^{1 - d/4}} \\
&= C' \sigma n^{-1/2} \alpha^{-d/8}.
\end{align*}
This completes the proof.
\end{proof}

With this estimate at hand, we now state the main convergence theorem for the fully discrete solution $f_h$ under sub-Gaussian
 noise.

\begin{theorem}\label{T5.2}
Assume Lemma~\ref{T3.1} holds, $\rho_0 = \|f^*\|_{L^2(\Omega)} + \sigma n^{-1/2}$, $\tau=O(h)$, $h=O(\alpha^{1/4})$, and $\alpha^{1/2+d/8} = O(\sigma n^{-1/2}\rho_0^{-1})$. Then for $z > 0$:
\begin{align}\label{4.22a}
\mathbb{P}\left( \|G_{\tau,h}f_h - Gf^*\|_n \geq \alpha^{1/2} \rho_0 z \right) &\leq 2e^{-C z^2}, \\\label{4.23a}
\mathbb{P}\left( \|f_h\|_{L^2(\Omega)} \geq \rho_0 z \right) &\leq 2e^{-C z^2}.
\end{align}
Moreover, if $\alpha \geq n^{-4/d}$, then
\begin{equation}\label{4.23aa}
\mathbb{P}\left( \|f_h - f^*\|_{H^{-1}(\Omega)} \geq \alpha^{1/4} \rho_0 z \right) \leq 2e^{-C z^2}.
\end{equation}
\end{theorem}
\begin{proof}
From a discrete Strang's lemma argument, we start with the inequality:
\[
\left\| \|f_n - f_h\|_{a_h} \right\|_{\psi_2} \leq C \left\| \inf_{v_h \in V_h} \|f_n - v_h\|_{a_h} \right\|_{\psi_2} + C \left\| \sup_{w_h \in V_h} \frac{|F(w_h)|}{\|w_h\|_{a_h}} \right\|_{\psi_2},
\]
where $\|v\|_{a_h}^2 $ and $F(w_h)$ are defined in \eqref{3.19a_full}.

We bound the two terms on the right. We bound the approximation error by choosing $v_h = P_h f_n$ and the estimate $\|f_n - P_h f_n\|_{a_h}^2 \le C(\alpha + h^4)\|f_n\|_{L^2(\Omega)}^2$. Taking the Orlicz norm and applying the assumption $h = O(\alpha^{1/4})$ along with the bound for $\|\|f_n\|\|_{\psi_2}$ from Theorem \ref{T4.1}, we obtain:

\[
\left\| \inf_{v_h \in V_h} \|f_n - v_h\|_{a_h} \right\|_{\psi_2} \le C\sqrt{\alpha+h^4} \|\|f_n\|_{L^2(\Omega)}\|_{\psi_2} \le C \alpha^{1/2}\rho_0.
\]
For the consistency error, we bound the three components of $F(w_h)$. The first two terms, which are deterministic and involve the operator difference $(G - G_{\tau,h})$, are bounded using Lemma \ref{L2.3a} and Theorem \ref{T4.1}. The assumption $h=O(\alpha^{1/4})$ ensures both are controlled by $C\alpha^{1/2}\rho_0$. From the proof of Theorem \ref{T5.1}, we know $\sup_{w_h \in V_h} |F(w_h)|/\|w_h\|_{a_h} = \sup_{w_h \in B_h} |F(w_h)|$, hence the random term is bounded by Lemma \ref{L4.6} and the crucial assumption on $\alpha$:
\[
\left\| \sup_{w_h} \frac{|(e, (G-G_{\tau,h}) w_h)_n|}{\|w_h\|_{a_h}} \right\|_{\psi_2} \le C\sigma n^{-1/2}\alpha^{-d/8} \le C(\rho_0\alpha^{1/2+d/8})\alpha^{-d/8} = C\rho_0\alpha^{1/2}.
\]
Combining these bounds, we conclude that the error between the discrete and continuous solutions is controlled:
\[
\left\| \|f_n - f_h\|_{a_h} \right\|_{\psi_2} \leq C\rho_0\alpha^{1/2}.
\]
The final results now follow. Using the triangle inequality, the bounds from Theorem \ref{T4.1}, and the result above, we can show $\|\|f_h\|_{L^2(\Omega)}\|_{\psi_2} \le C\rho_0$ and $\|\| G_{\tau,h}f_h - Gf^* \|_n \|_{\psi_2} \le C\alpha^{1/2}\rho_0$. The $H^{-1}$ estimate \eqref{4.23aa} follows from the same line of reasoning as in Theorem \ref{T4.2}, now applied to the discrete solution $f_h$. Applying the property of the $\psi_2$-norm to these Orlicz norm bounds gives the stated probability estimates.
\end{proof}

\section{Numerical experiment}\label{4}

 Motivated by Theorems~\ref{T5.1} and~\ref{T5.2}, we introduce a self-consistent algorithm to determine the regularization parameter $\alpha$ in problem~\eqref{5.6}. The algorithm is derived from the heuristic rule
\begin{equation}\label{6.1}
\alpha^{1/2 + d/8} = \sigma n^{-1/2} \|f^*\|_{L^2(\Omega)}^{-1},
\end{equation}
which is designed to balance the bias--variance trade-off in the reconstruction error while ensuring consistency in the $H^{-1}$-norm between the exact and reconstructed solutions.  In this section, we denote the relative noise level as $\delta = 
\sigma / \|G f^{*}\|_{L^{\infty}(\Omega)}
$.

In practice, however, the true source $f^*$ and noise level $\sigma$ are generally unknown. To address this, we propose an iterative scheme that self-consistently approximates $\|f^*\|_{L^2(\Omega)}$ using the computed solution $\|f_h\|_{L^2(\Omega)}$, and estimates the noise level $\sigma$ from the residual norm $\|G_{\tau,h}f_h - m\|_n$; this estimation relies on the identity $\|Gf^* - m\|_n = \|e\|_n$. For the numerical validation, the $H^{-1}(\Omega)$-norm error $\|f^* - f_h\|_{H^{-1}(\Omega)}$ is evaluated using the technique described in \cite[Section 6]{38}, where it is approximated by $\|P_h f^* - f_h\|_{H^{-1}(\Omega)}$.

\begin{Algo}
[Estimation of the regularization parameter $\alpha$.]
\begin{enumerate}
    \item Initialize: Set initial guess $\alpha_{0}= n^{-4/(d+4)}$.
    \item While $\left| \frac{\alpha_{j-1} - \alpha_{j}}{\alpha_{j}} \right| > 10^{-3}$:
    \begin{enumerate}
        \item Solve~\eqref{5.6} for $f_h$ using $\alpha = \alpha_{j}$ on mesh $\mathcal{M}_h$.
        \item Update the regularization parameter:
        \[
        \alpha_{j+1}^{1/2 + d/8} = n^{-1/2} \|G_{\tau,h}f_h - m\|_n \cdot \|f_h\|_{L^2(\Omega)}^{-1}.
        \]
    \end{enumerate}
\end{enumerate}
\label{Algo1}
\end{Algo}

\begin{example}
Consider the wave equation with \( f(x) = x^{1/4}(1 - x)^{1/4} \), $g(t)=t^4$ and \( \Omega = (0,1) \).
\end{example}

\begin{figure}[htbp]
\centering
\subfigure[$\|Gf^* - G_{\tau,h}f_h\|_n$]{
\includegraphics[width=4.5cm]{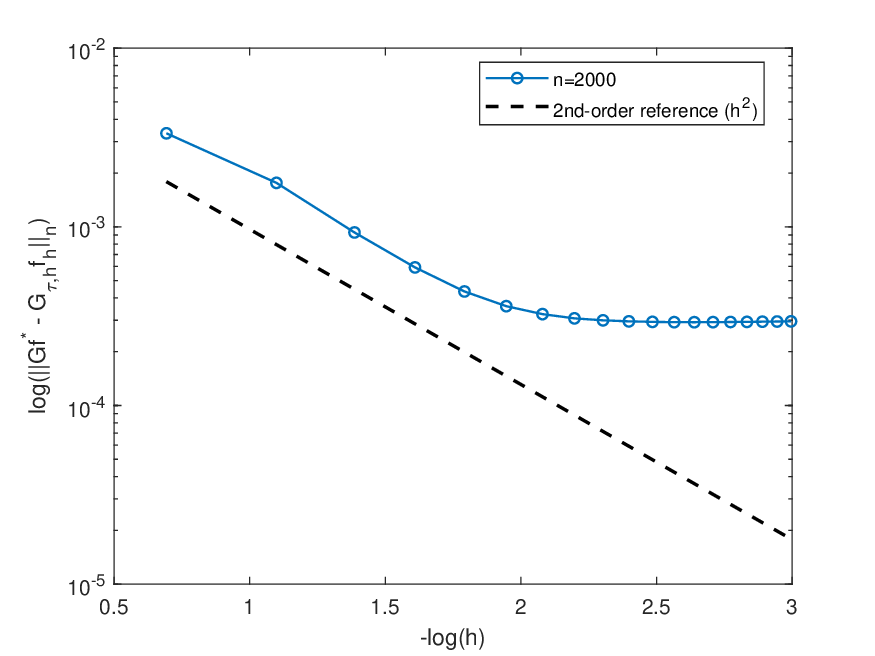}
}
\subfigure[ $\|f^* - f_h\|_{H^{-1}(\Omega)}$]{
\includegraphics[width=4.5cm]{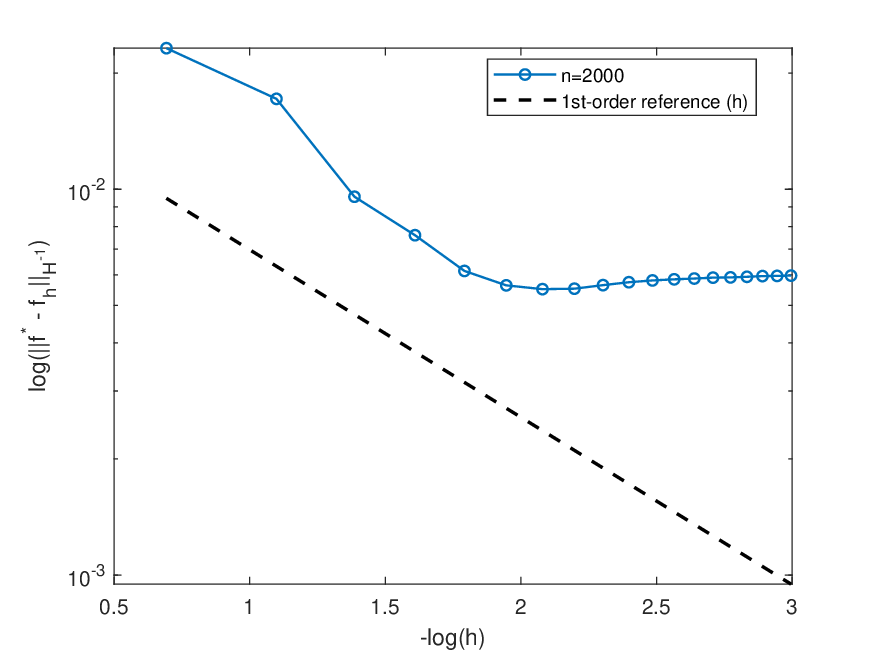}
}
\caption{
The reference lines indicate the expected $h^2$ and $h$ convergence rates for the empirical norm and $H^{-1}$ norm errors, respectively.
}\label{F41a}
\end{figure}

We begin with the one-dimensional case ($d=1$). The convergence of the finite element method is examined for a fixed $\alpha = 10^{-6}$, with $\sigma$ and $n$ chosen following \eqref{6.1}. Figure~\ref{F41a} displays the resulting convergence rates: approximately $\mathcal{O}(h^{2})$ for the empirical norm error and $\mathcal{O}(h^{1})$ for the $H^{-1}$ norm error. These rates are observed until the modeling error, governed by $\alpha$, becomes dominant. A good agreement is observed between these numerical results and the theoretical estimates \eqref{3.13ab} and \eqref{step3a}.

\begin{figure}[htbp]
\centering
\subfigure[$\|Gf^* - G_{\tau,h}f_h\|_n$]{
\includegraphics[width=4.5cm]{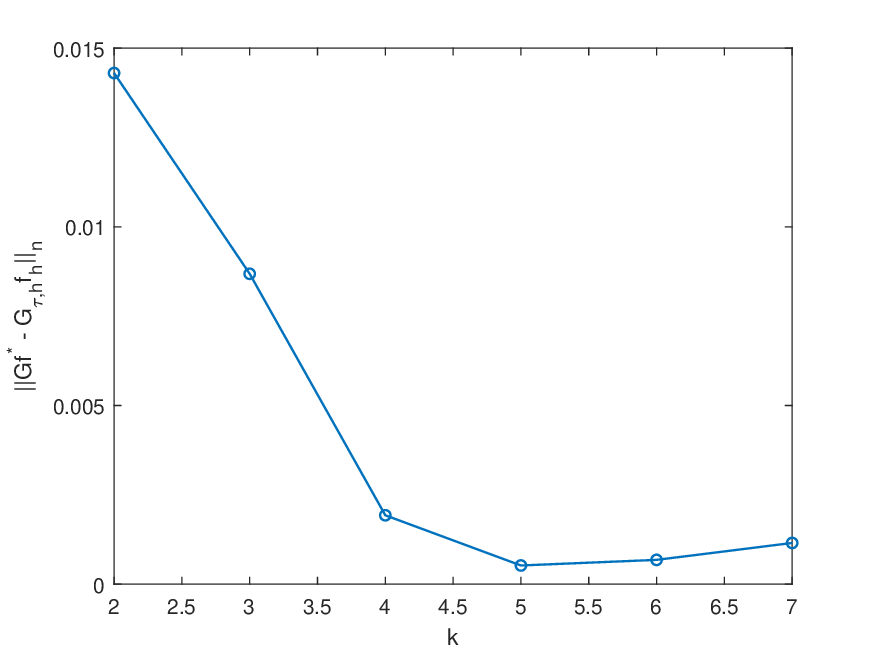}
}
\subfigure[ $\|f^* - f_h\|_{H^{-1}(\Omega)}$]{
\includegraphics[width=4.5cm]{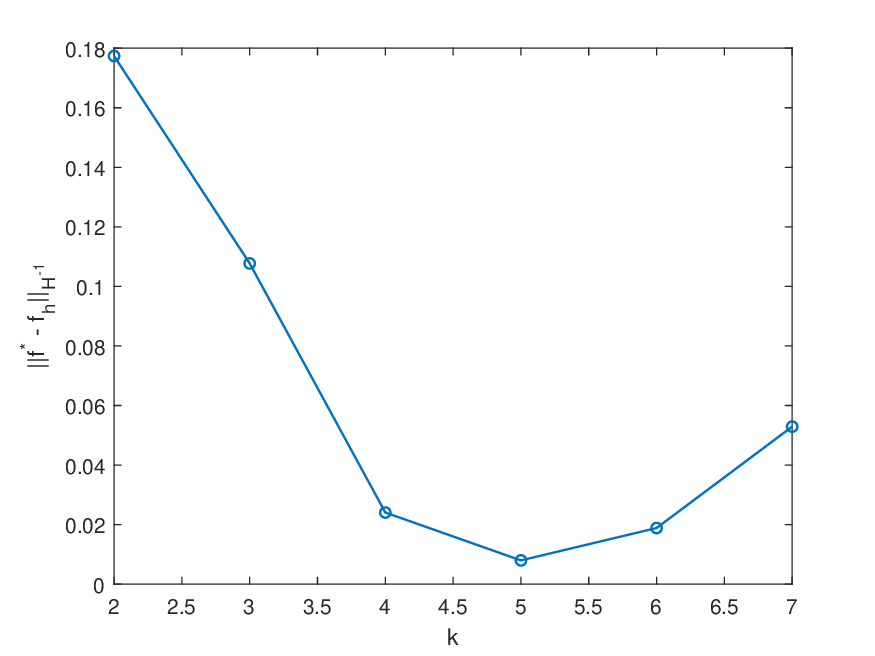}
}
\caption{
Dependence of (a) the empirical error and (b) the $H^{-1}(\Omega)$-norm error on the regularization parameter $\alpha = 10^{-k}$ ($k = 2, \ldots, 7$), for the case with a noise level of $\delta = 40.81\%$ and $n = 300$ 
measurements.
}
\label{F4.2o}
\end{figure}

To verify the near-optimality of the smoothing parameter selection rule for $\alpha$ given in \eqref{6.1}, we choose the computational parameters as follows: mesh size $h = 1/251$, and time step $\tau = 1/200$. These values are sufficiently small to ensure that discretization errors from the finite element approximation are negligible. For a noise level $\sigma = 0.009$ with $\|G f^{*}\|_{L^{\infty}(\Omega)}  \approx 0.0221$ (yielding relative noise level $\delta =40.81\%$) and  sample size $n = 300$, the rule in \eqref{6.1} gives $\alpha \approx 1.1749 \times 10^{-5}$. This value is consistent with the empirical optimum $\alpha = 1 \times 10^{-5}$ obtained by minimizing the error norms over $\alpha = 10^{-k}$ ($k = 2, \dots, 7$), as shown in Figures~\ref{F4.2o}(a)--(b).


\begin{figure}[htbp]
\centering
\subfigure[$\|Gf^* - G_{\tau,h}f_h\|_n$]{
\includegraphics[width=4.5cm]{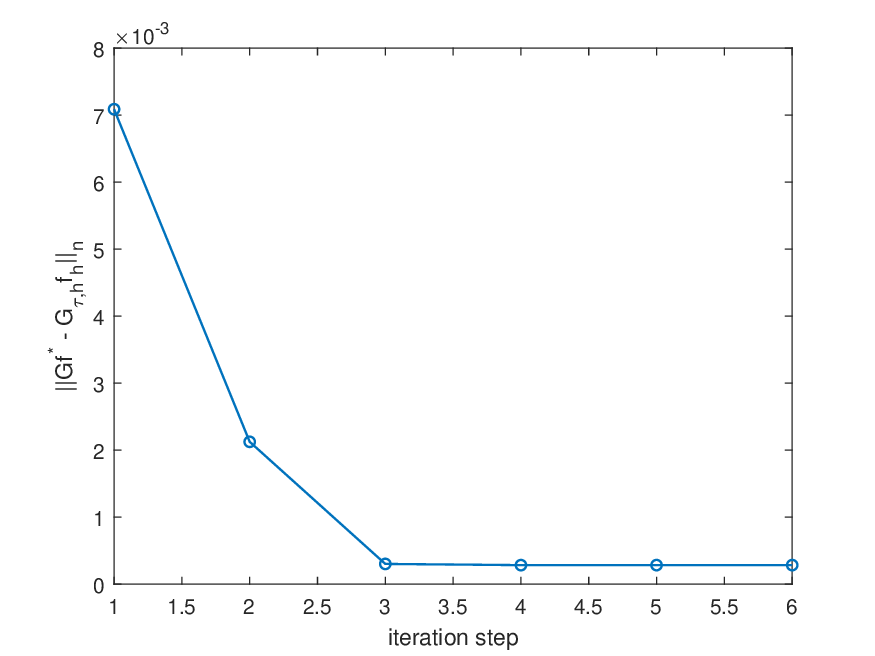}
}
\subfigure[$\|f^* - f_h\|_{H^{-1}(\Omega)}$]{
\includegraphics[width=4.5cm]{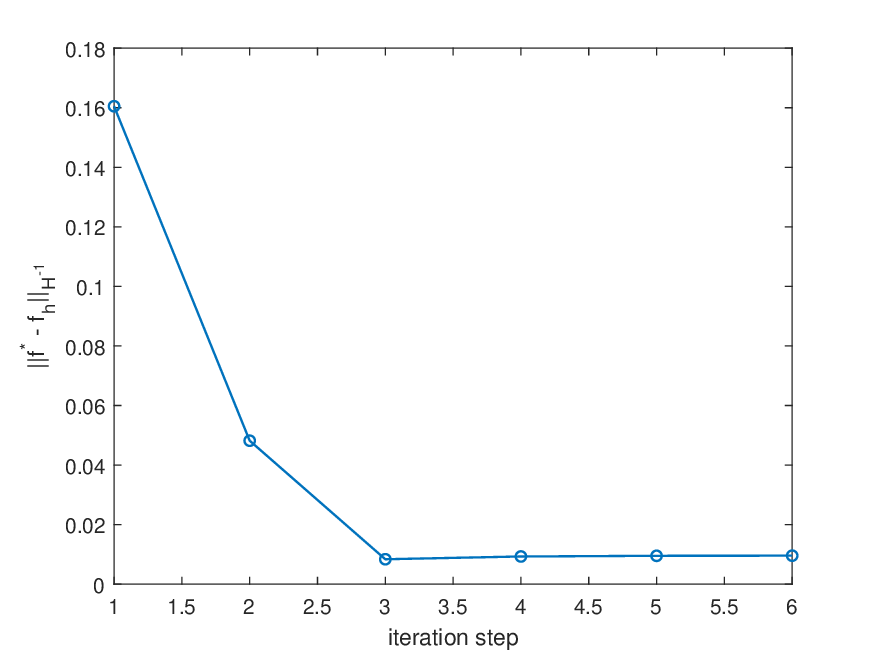}
}
  \caption{
Evolution of the empirical error $\|Gf^* - G_{\tau,h}f_h\|_n$ and $\|f^* - f_h\|_{H^{-1}(\Omega)}$ across iterations for a fixed noise level $\delta = 0.4081$ with  $n = 1000$ measurement points.
}
  \label{oF4.4}
\end{figure}

Subsequently, Algorithm~\ref{Algo1} is employed to compute the regularization parameter $\alpha$. For the case with $\sigma = 0.009$ and $n = 1000$, Figure~\ref{oF4.4}(b) illustrates the convergence behavior of the sequence $\{\alpha_{k}\}$. The computed value $\alpha_{6} = 4.2989 \times 10^{-6}$ agrees well with the theoretical optimum $4.4845 \times 10^{-6}$ given by \eqref{6.1}. Moreover, the residual estimate $\|m - G_{\tau,h} f_h\|_n = 0.0091$ 
is consistent with the specified noise level.

\begin{figure}[htbp]
\centering
\includegraphics[width=4.5cm]{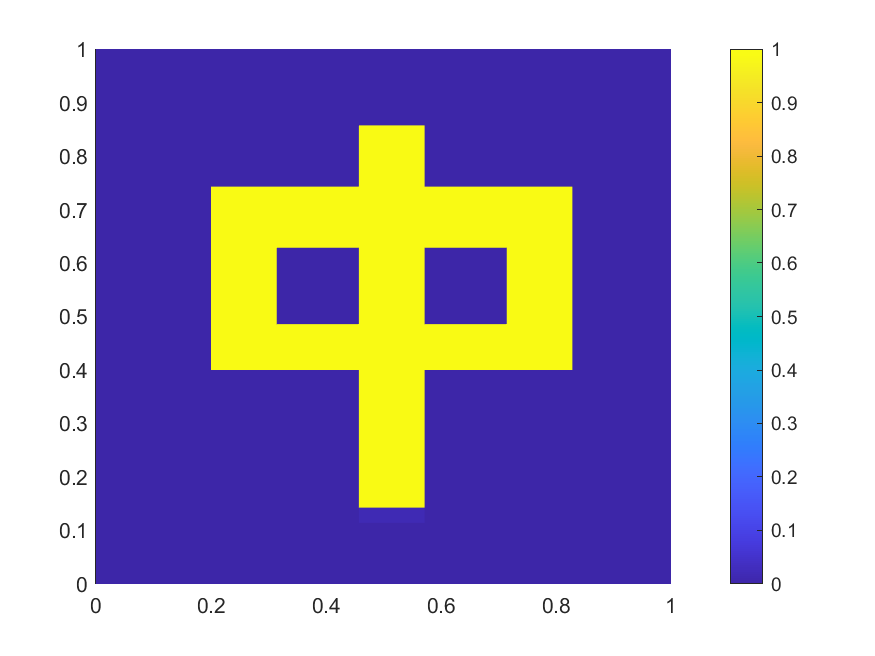}
\caption{Surface plot of the exact source term $f^*$.}\label{F1}
\end{figure}

\begin{example}\label{exam1}
Consider the wave equation defined on the domain $\Omega = (0,1) \times (0,1)$, with the source function $f$ illustrated in Figure~\ref{F1}.
\end{example}

\begin{figure}[htbp]
\centering
\subfigure[$\|Gf^* - G_{\tau,h}f_h\|_n$]{
\includegraphics[width=4.5cm]{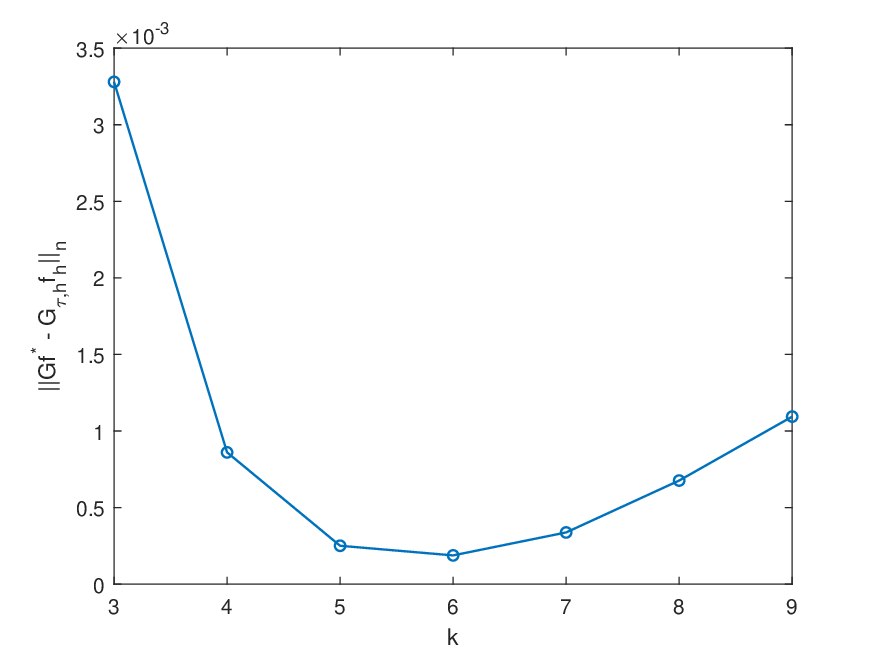}
}
\subfigure[$\|f^* - f_h\|_{H^{-1}(\Omega)}$]{
\includegraphics[width=4.5cm]{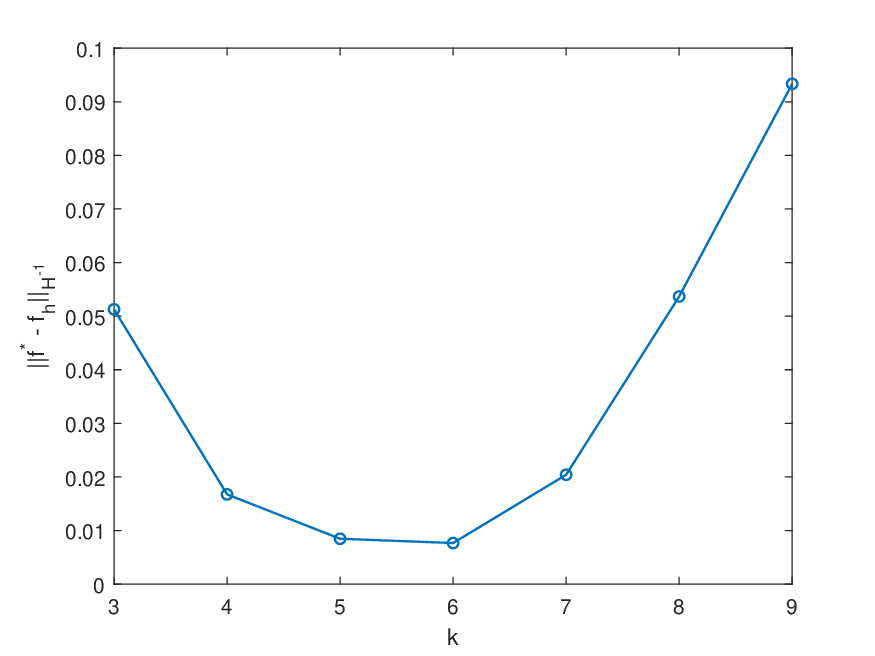}
}
\caption{Dependence of (a) the empirical error and (b) the $H^{-1}(\Omega)$-norm error on the regularization parameter $\alpha = 10^{-k}$ ($k = 3, \ldots, 9$), for the case with a noise level of $\delta = 11.1\%$ and $n = 2500$ measurements.}\label{F4.2}
\end{figure}

For the two-dimensional case ($d=2$), we first verify the near-optimality of the smoothing parameter selection rule in \eqref{6.1}. We set $n = 2500$, mesh size $h = 1/31$, and time step $\tau = 1/200$ to ensure that finite element discretization errors are negligible. The maximum value of the exact solution is $\|G f^{*}\|_{L^{\infty}(\Omega)}  \approx 0.018$. Given $\|f^*\|_{L^2(\Omega)} \approx 0.4614$ and noise level $\sigma = 0.002$ (corresponding to $\delta = 11.1\%$), the rule yields $\alpha \approx 3.837 \times 10^{-6}$, which is consistent with the empirically optimal value $\alpha = 1.0 \times 10^{-6}$; see Figures~\ref{F4.2}(a)--(b).

\begin{figure}[htbp]
\centering
\subfigure[\( \|G f^* - G_{\tau, h} f_h\|_n \) . ]{
\includegraphics[width=4.5cm]{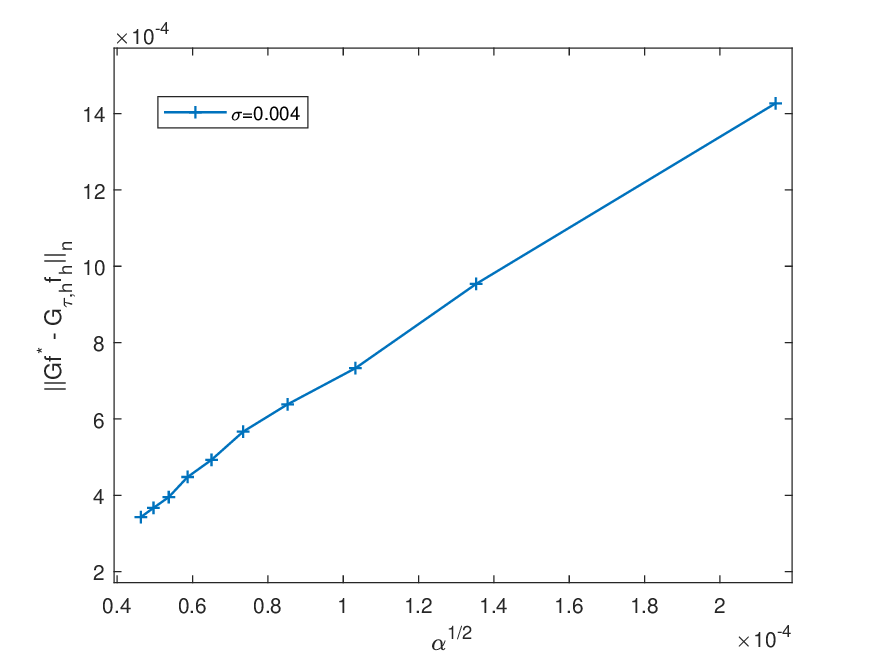}
}
\subfigure[\( \|f^* - f_h\|_{H^{-1}(\Omega)} \) .  ]{
\includegraphics[width=4.5cm]{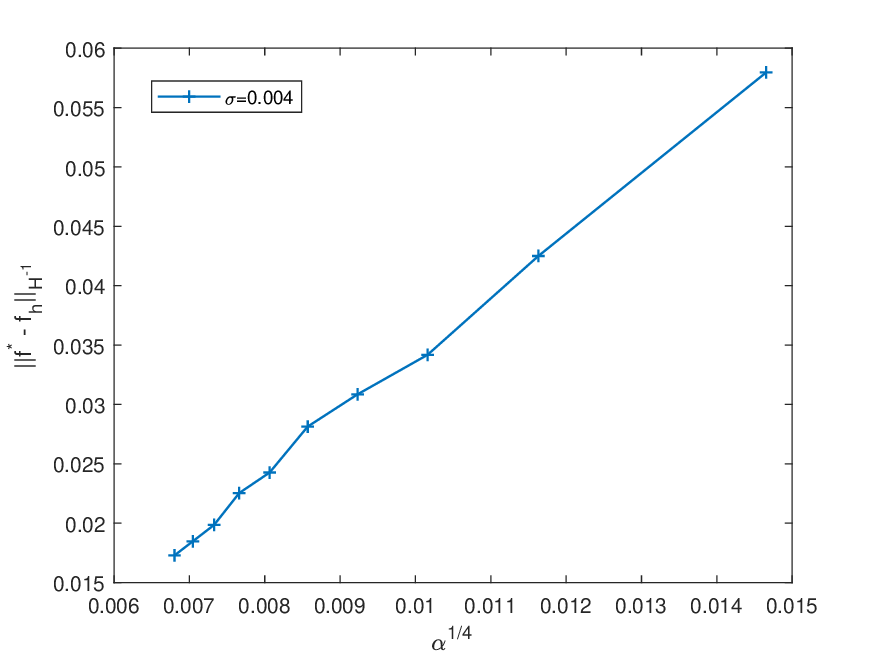}
}
\caption{Linear dependence of the empirical approximation error on $\alpha^{1/2}$ (a) and of the $H^{-1}(\Omega)$-norm error on $\alpha^{1/4}$ (b), for the case with $\sigma = 0.004$ and a relative noise level of $\delta = 22.2\%$.}\label{F4.3}
\end{figure}

To further validate Theorems~\ref{T5.1} and \ref{T5.2}, we investigate whether the empirical approximation error $\|G f^* - G_{\tau, h} f_h\|_n$ exhibits linear dependence on $\alpha^{1/2}$ and the $H^{-1}(\Omega)$-norm error $\|f^* - f_h\|_{H^{-1}(\Omega)}$ shows linear dependence on $\alpha^{1/4}$ when $\alpha$ is selected via the optimal rule \eqref{6.1}. The true source $f^*$ is defined in Figure~\ref{F1}, with $n$ ranging from $25 \times 10^2$ to $25 \times 10^4$. For $\sigma = 0.004$, the relative noise level is $\delta = 22.2\%$. Figure~\ref{F4.3} demonstrates the linear dependence of $\|G f^* - G_{\tau, h} f_h\|_n$ on $\alpha^{1/2}$ and $\|f^* - f_h\|_{H^{-1}(\Omega)}$ on $\alpha^{1/4}$.

\begin{figure}[htbp]
\centering
\subfigure[$\|G f^* - G_{\tau, h} f_h\|_n$ convergence]{
\includegraphics[width=4.5cm]{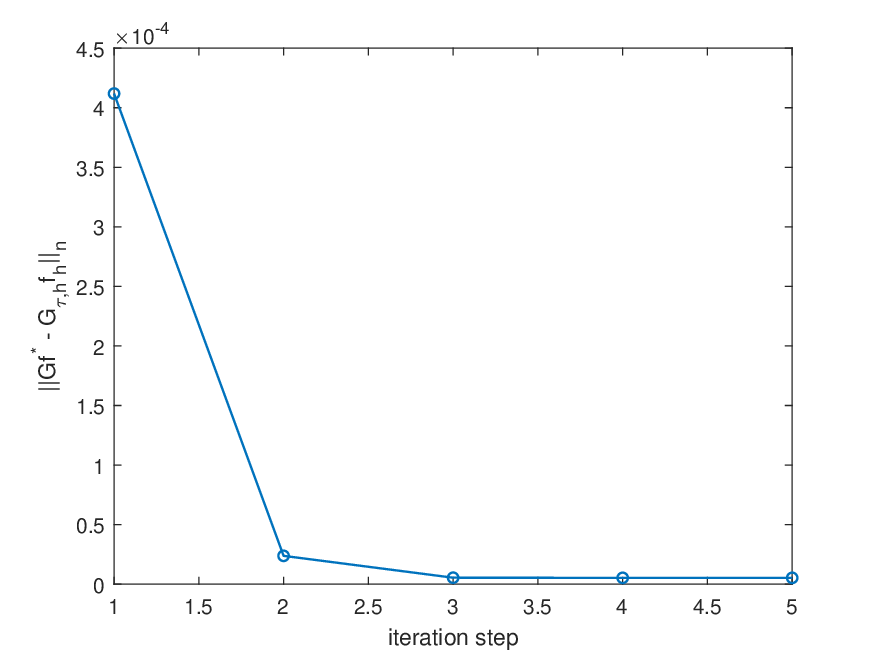}
}
\subfigure[Iteration  1]{
\includegraphics[width=4.5cm]{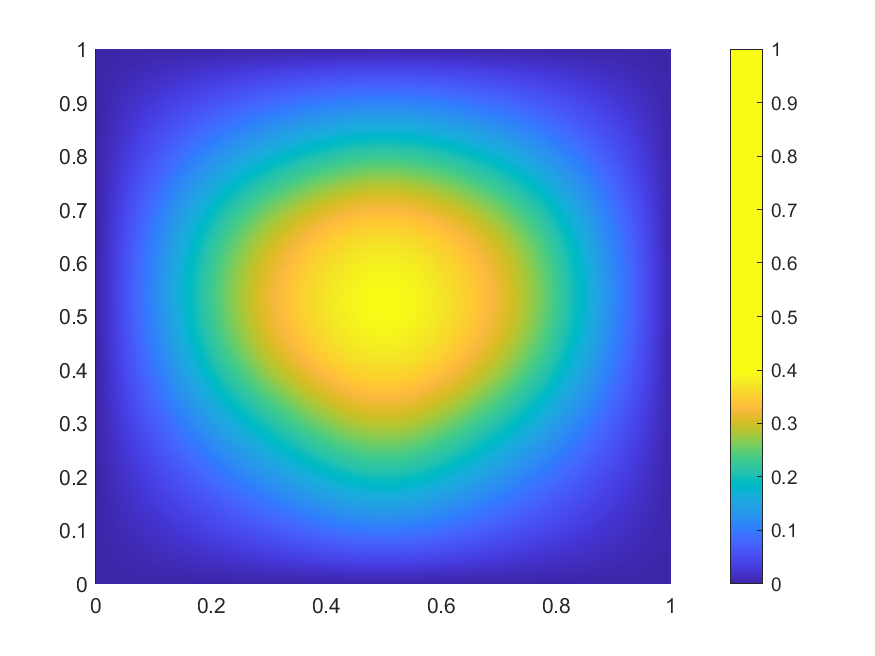}
}
\subfigure[Iteration 5]{
\includegraphics[width=4.5cm]{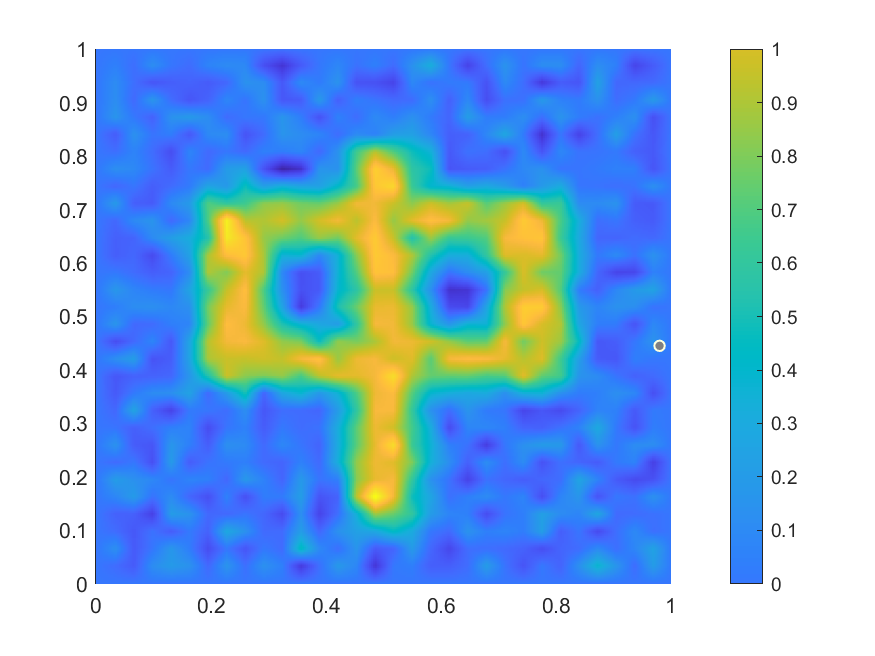}
}
\caption{Results for fixed $\sigma = 0.001$, $\delta = 5.55\%$, and $n=90000$: (a) Evolution of the empirical error $\|Gf^* - G_{\tau,h}f_h\|_n$ with iterations; (b)--(c) Reconstructions of $f_h$ at iterations 1 and 5.} \label{F4.4}
\end{figure}

Finally, having validated the near-optimality of the selection rule \eqref{6.1} for $\alpha$ and the  convergence (Theorem~\ref{T5.2}), we evaluate the performance of Algorithm~\ref{Algo1} in determining the regularization parameter $\alpha$ without prior knowledge of $f^*$ or $\sigma$. We set $n = 9 \times 10^4$ and $\sigma = 0.001$. Figure~\ref{F4.4}(a) shows the convergence of $\{\alpha_{k}\}$, where the computed value $\alpha_{5} = 1.4344 \times 10^{-7}$ aligns well with the theoretical optimum $1.3966 \times 10^{-7}$ from \eqref{6.1}. The reconstructed source term $f_h$ at the fifth iteration (Figure~\ref{F4.4}(c)) shows significantly improved accuracy and better approximates the true solution compared to the first iteration result (Figure~\ref{F4.4}(b)). Moreover, the residual estimate $\|m - G_{\tau,h} f_h\|_n \approx 9.97\times 10^{-4}$ matches the prescribed noise level $\sigma$.

\vspace{3em}

\section*{Acknowledgments.}
Zhidong Zhang is supported by the National Key Research and Development Plan of China (Grant No.2023YFB3002400). Wenlong Zhang is partially supported by the National Natural Science Foundation of China under grant numbers No.12371423 and No.12241104.

\end{document}